\documentclass{article}
\usepackage{amsmath}
\usepackage{amssymb}

\usepackage{theorem}
\numberwithin{equation}{section}

\newtheorem{theorem}{Theorem}[section]

\newtheorem{corollary}[theorem]{Corollary}
\newtheorem{lemma}[theorem]{Lemma}
\newtheorem{definition}{Definition}[section]

\begin{document}
\title{Global well-posedness and scattering for the radial, defocusing, cubic nonlinear wave equation}
\date{\today}
\author{Benjamin Dodson}
\maketitle

\noindent \textbf{Abstract:} In this paper we prove global well-posedness and scattering for the defocusing, cubic, nonlinear wave equation on $\mathbf{R}^{1 + 3}$ with radial initial data lying in the critical Sobolev space $\dot{H}^{1/2}(\mathbf{R}^{3}) \times \dot{H}^{-1/2}(\mathbf{R}^{3})$.

\section{Introduction}
In this paper we study the defocusing, cubic nonlinear wave equation
\begin{equation}\label{1.1}
u_{tt} - \Delta u + u^{3} = 0, \qquad u(0,x) = u_{0}, \qquad u_{t}(0,x) = u_{1}.
\end{equation}
This problem is $\dot{H}^{1/2}$ critical, since the equation $(\ref{1.1})$ is invariant under the scaling symmetry
\begin{equation}\label{1.2}
u(t,x) \mapsto \lambda u(\lambda t, \lambda x).
\end{equation}
This scaling symmetry completely determines local well-posedness theory for $(\ref{1.1})$. Positively, \cite{LS} proved
\begin{theorem}\label{t1.1}
The equation $(\ref{1.1})$ is locally well-posed for initial data in $u_{0} \in \dot{H}^{1/2}(\mathbf{R}^{3})$ and $u_{1} \in \dot{H}^{-1/2}(\mathbf{R}^{3})$ on some interval $[-T(u_{0}, u_{1}), T(u_{0}, u_{1})]$. The time of well-posedness $T(u_{0}, u_{1})$ depends on the profile of the initial data $(u_{0}, u_{1})$, not just its size. 

Additional regularity is enough to give a lower bound on the time of well-posedness. Therefore, there exists some $T(\| u_{0} \|_{\dot{H}^{s}}, \| u_{1} \|_{\dot{H}^{s - 1}}) > 0$ for any $\frac{1}{2} < s < \frac{3}{2}$.
\end{theorem}
Negatively, \cite{LS} proved
\begin{theorem}\label{t1.1.1}
Equation $(\ref{1.1})$ is ill-posed for $u_{0} \in \dot{H}^{s}(\mathbf{R}^{3})$ and $u_{1} \in \dot{H}^{s - 1}(\mathbf{R}^{3})$ when $s < \frac{1}{2}$.
\end{theorem}
Local well-posedness is defined in the usual way.
\begin{definition}[Locally well-posed]\label{d1.2}
The initial value problem $(\ref{1.1})$ is said to be locally well-posed if there exists an open interval $I \subset \mathbf{R}$ containing $0$ such that:
\begin{enumerate}
\item A unique solution $u \in L_{t}^{\infty} \dot{H}^{1/2}(I \times \mathbf{R}^{3}) \cap L_{t,loc}^{4} L_{x}^{4}(I \times \mathbf{R}^{3})$, $u_{t} \in L_{t}^{\infty} \dot{H}^{-1/2}(I \times \mathbf{R}^{3})$ exists.

\item The solution $u$ is continuous in time, $u \in C(I ; \dot{H}^{1/2}(\mathbf{R}^{3}))$, $u_{t} \in C(I ; \dot{H}^{-1/2}(\mathbf{R}^{3}))$.

\item The solution $u$ depends continuously on the initial data in the topology of item one.
\end{enumerate}
\end{definition}

Given this fact, it is natural to inquire as to the long-time behavior of solutions to $(\ref{1.1})$ with initial data at the $\dot{H}^{1/2}$-critical regularity. Do they continue for all time, and if they do, what is their behavior at large times?

Global well-posedness for initial data in $\dot{H}^{1/2} \cap \dot{H}^{1}(\mathbf{R}^{3}) \times \dot{H}^{-1/2} \cap L^{2}(\mathbf{R}^{3})$ follows from conservation of the energy
\begin{equation}\label{1.4}
E(u(t)) = \frac{1}{2} \int u_{t}(t,x)^{2} dx + \frac{1}{2} \int |\nabla u(t,x)|^{2} dx + \frac{1}{4} \int u(t,x)^{4} dx.
\end{equation}
By the Sobolev embedding theorem and H\"older's inequality,
\begin{equation}\label{1.4.1}
\| u(0) \|_{L_{x}^{4}(\mathbf{R}^{3})}^{4} \lesssim \| u(0) \|_{L_{x}^{3}(\mathbf{R}^{3})}^{2} \| u(0) \|_{L_{x}^{6}(\mathbf{R}^{3})}^{2} \lesssim \| u(0) \|_{\dot{H}^{1/2}(\mathbf{R}^{3})}^{2} \| u(0) \|_{\dot{H}^{1}(\mathbf{R}^{3})}^{2},
\end{equation}
and therefore,
\begin{equation}\label{1.4.2}
E(u(0)) \lesssim_{\| u_{0} \|_{\dot{H}^{1/2}}} \| u_{0} \|_{\dot{H}^{1}(\mathbf{R}^{3})}^{2} + \| u_{1} \|_{L^{2}(\mathbf{R}^{3})}^{2}.
\end{equation}
By $(\ref{1.4})$, $E(u(t)) = E(u(0))$ controls the size of $\| u(t) \|_{\dot{H}^{1}} + \| u_{t}(t) \|_{L^{2}}$, which by Theorem $\ref{t1.1}$ gives global well-posedness.

Comparing $(\ref{1.1})$ to the quintic wave equation in three dimensions,
\begin{equation}\label{1.4.3}
u_{tt} - \Delta u + u^{5} = 0, \qquad u(0,x) = u_{0}, \qquad u_{t}(0,x) = u_{1},
\end{equation}
a solution to $(\ref{1.4.3})$ is invariant under the scaling symmetry $u(t,x) \mapsto \lambda^{1/2} u(\lambda t, \lambda x)$, a symmetry that preserves the $\dot{H}^{1} \times L^{2}$ norm of $(u_{0}, u_{1})$. Observe that the conserved energy for $(\ref{1.4.3})$,
\begin{equation}\label{1.4.4}
E(u(t)) = \frac{1}{2} \int u_{t}(t,x)^{2} dx + \frac{1}{2} \int |\nabla u(t,x)|^{2} dx + \frac{1}{6} \int u(t,x)^{6} dx
\end{equation}
is also invariant under the scaling symmetry. For this reason, $(\ref{1.4.3})$ is called energy--critical, and it is possible to prove a result in the same vein as Theorems $\ref{t1.1}$ and $\ref{t1.1.1}$ at the critical regularity $\dot{H}^{1} \times L^{2}$. 

This fact combined with conservation of the energy $(\ref{1.4.4})$ is insufficient to prove global well-posedness for $(\ref{1.4.3})$. The reason is because the time of local well-posedness depends on the profile of the initial data $(u_{0}, u_{1}) \in \dot{H}^{1} \times L^{2}$, and not just its size. Instead, the proof of global well-posedness for the quintic problem uses a non-concentration of energy argument. This result has been completely worked out, proving both global well-posedness and scattering, for both the radial (\cite{GSV}, \cite{Struwe}) and the nonradial case (\cite{BS}, \cite{Gril}, \cite{Shatah - Struwe}).
\begin{definition}[Scattering]\label{d1.1}
A solution to $(\ref{1.4.3})$ is said to be scattering in some $\dot{H}^{s}(\mathbf{R}^{3}) \times \dot{H}^{s - 1}(\mathbf{R}^{3})$ if there exist $(u_{0}^{+}, u_{1}^{+}), (u_{0}^{-}, u_{1}^{-}) \in \dot{H}^{s} \times \dot{H}^{s - 1}$ such that
\begin{equation}\label{1.4.5}
\lim_{t \rightarrow +\infty} \| (u(t), u_{t}(t)) - S(t)(u_{0}^{+}, u_{1}^{+}) \|_{\dot{H}^{s} \times \dot{H}^{s - 1}} = 0,
\end{equation}
and
\begin{equation}\label{1.4.6}
\lim_{t \rightarrow -\infty} \| (u(t), u_{t}(t)) - S(t)(u_{0}^{+}, u_{1}^{+}) \|_{\dot{H}^{s} \times \dot{H}^{s - 1}} = 0,
\end{equation}
where $S(t)(f, g)$ is the solution operator to the linear wave equation. That is, if $(u(t), u_{t}(t)) = S(t)(f, g)$, then
\begin{equation}\label{1.4.7}
u_{tt} - \Delta u = 0, \qquad u(0,x) = f, \qquad u_{t}(0,x) = g.
\end{equation}
\end{definition}

Similar results for $(\ref{1.1})$ may also be obtained if one assumes a uniform bound over $\| u \|_{\dot{H}^{1/2}(\mathbf{R}^{3})} + \| u_{t} \|_{\dot{H}^{-1/2}(\mathbf{R}^{3})}$ for the entire time of existence of the solution.
\begin{theorem}\label{t1.3}
Suppose $u_{0} \in \dot{H}^{1/2}(\mathbf{R}^{3})$ and $u_{1} \in \dot{H}^{-1/2}(\mathbf{R}^{3})$ are radial functions, and $u$ solves $(\ref{1.1})$ on a maximal interval $0 \in I \subset \mathbf{R}$, with
\begin{equation}\label{1.6}
\sup_{t \in I} \| u(t) \|_{\dot{H}^{1/2}(\mathbf{R}^{3})} + \| u_{t}(t) \|_{\dot{H}^{-1/2}(\mathbf{R}^{3})} < \infty.
\end{equation}
Then $I = \mathbf{R}$ and the solution $u$ scatters both forward and backward in time.
\end{theorem}
\emph{Proof:} See \cite{DL}. $\Box$\vspace{5mm}

In this paper we remove the a priori assumption on uniform boundedness of the critical norm in $(\ref{1.6})$, proving,
\begin{theorem}\label{t1.2}
The initial value problem $(\ref{1.1})$ is globally well-posed and scattering for radial initial data $u_{0} \in \dot{H}^{1/2}(\mathbf{R}^{3})$ and $u_{1} \in \dot{H}^{-1/2}(\mathbf{R}^{3})$. Moreover, there exists a function $f : [0, \infty) \rightarrow [0, \infty)$ such that if $u$ solves $(\ref{1.1})$ with initial data $(u_{0}, u_{1}) \in \dot{H}^{1/2} \times \dot{H}^{-1/2}$, then
\begin{equation}\label{1.5}
\| u \|_{L_{t,x}^{4}(\mathbf{R} \times \mathbf{R}^{3})} \leq f(\| u_{0} \|_{\dot{H}^{1/2}(\mathbf{R}^{3})} + \| u_{1} \|_{\dot{H}^{-1/2}(\mathbf{R}^{3})}).
\end{equation}
\end{theorem}

The proof of Theorem $\ref{t1.2}$ combines the Fourier truncation method and hyperbolic coordinates. Previously, \cite{KPV} applied the Fourier truncation method to the cubic wave equation, $(\ref{1.1})$, proving global well-posedness of $(\ref{1.1})$ with initial data lying in the inhomogeneous Sobolev spaces $H_{x}^{s}(\mathbf{R}^{3}) \times H_{x}^{s - 1}(\mathbf{R}^{3})$ for $s > \frac{3}{4}$. This argument was improved and modified in many subsequent papers, for both radial and nonradial data. In particular, see \cite{D} for a proof of global well-posedness for $(\ref{1.1})$ with radial initial data lying in
\begin{equation}\label{1.7}
(\dot{H}^{s}(\mathbf{R}^{3}) \cap \dot{H}^{1/2}(\mathbf{R}^{3})) \times (\dot{H}^{s - 1}(\mathbf{R}^{3}) \cap \dot{H}^{-1/2}(\mathbf{R}^{3})),
\end{equation}
for any $s > \frac{1}{2}$, as well as for a description of other results along this line.\vspace{5mm}

\noindent \textbf{Remark:} The method used in \cite{D} was the I-method, a modification of the Fourier truncation method.\vspace{5mm}

In this paper, using the Fourier truncation method, global well-posedness is proved for $(\ref{1.1})$ with radial initial data lying in $\dot{H}^{1/2}(\mathbf{R}^{3}) \times \dot{H}^{-1/2}(\mathbf{R}^{3})$. The idea behind the proof is that at low frequencies, the initial data has finite energy, and a solution to $(\ref{1.1})$ with finite energy is global. Meanwhile, at high frequencies, the $\dot{H}^{1/2} \times \dot{H}^{-1/2}$ norm is small, and for such initial data, $(\ref{1.1})$ may be treated using perturbative arguments. The mixed terms in the nonlinearity are then shown to have finite energy, proving global well-posedness.

Proof of scattering utilizes hyperbolic coordinates. Hyperbolic coordinates were used in \cite{Tataru} to prove weighted Strichartz estimates that were proved in \cite{GLS}. More recently, \cite{Shen}, working in hyperbolic coordinates, was able to prove a scattering result for data lying in a weighted energy space. Later, \cite{D1} combined the result of \cite{Shen} with the I-method argument in \cite{D} to prove scattering data lying in the subspace of $\dot{H}^{1/2} \times \dot{H}^{-1/2}$,
\begin{equation}\label{1.8}
\| u_{0} \|_{\dot{H}^{1/2 + \epsilon}(\mathbf{R}^{3})} + \| |x|^{2 \epsilon} u_{0} \|_{\dot{H}^{1/2 + \epsilon}(\mathbf{R}^{3})} + \| u_{1} \|_{\dot{H}^{-1/2 + \epsilon}(\mathbf{R}^{3})} + \| |x|^{2 \epsilon} u_{1} \|_{\dot{H}^{-1/2 + \epsilon}(\mathbf{R}^{3})}.
\end{equation}
Here, the Fourier truncation global well-posedness argument in hyperbolic coordinates shows that $(\ref{1.1})$ is globally well-posed and scattering for any $(u_{0}, u_{1}) \in \dot{H}^{1/2} \times \dot{H}^{-1/2}$.

This fact still falls short of $(\ref{1.5})$, since the proof does not give any uniform control over the $\| u \|_{L_{t,x}^{4}(\mathbf{R} \times \mathbf{R}^{3})}$ norm. To remedy this deficiency, and complete the proof of Theorem $\ref{t1.2}$, a profile decomposition is used. The profile decomposition shows that for any bounded sequence of initial data
\begin{equation}\label{1.9}
\| u_{0}^{n} \|_{\dot{H}^{1/2}(\mathbf{R}^{3})} + \| u_{1}^{n} \|_{\dot{H}^{-1/2}(\mathbf{R}^{3})} \leq A,
\end{equation}
and if $u^{n}(t)$ is the global solution to $(\ref{1.1})$ with initial data $(u_{0}^{n}, u_{1}^{n})$, then
\begin{equation}\label{1.10}
\| u^{n} \|_{L_{t,x}^{4}(\mathbf{R} \times \mathbf{R}^{3})} < \infty,
\end{equation}
is uniformly bounded. Then by Zorn's lemma, the proof of Theorem $\ref{t1.2}$ is complete.\vspace{5mm}

The author believes this to be the first unconditional global well-posedness and scattering result for a nonlinear wave equation with initial data lying in the critical Sobolev space, with no conserved quantity that controls the critical norm. Previously, \cite{D2} proved global well-posedness and scattering for $(\ref{1.1})$ with radial initial data lying in the Besov space $B_{1,1}^{2} \times B_{1,1}^{1}$. These spaces are also invariant under the scaling $(\ref{1.2})$. Later, \cite{MYZ} proved a similar result in five dimensions.

There are two main improvements for this result over the results of \cite{D2} and \cite{MYZ}. The first is that, while scale invariant, the Besov spaces are only subsets of the critical Sobolev spaces. The second improvement is that the $\dot{H}^{1/2} \times \dot{H}^{-1/2}$ norm is invariant under the free evolution of the linear wave equation. Whereas, for initial data lying in a Besov space, the proof of scattering simply meant that the solution scattered in the $\dot{H}^{1/2} \times \dot{H}^{-1/2}$ norm.\vspace{5mm}

\noindent \textbf{Acknowledgements:} The author was partially supported on NSF grant number $1764358$ during the writing of this paper. The author was also a guest of the Institute for Advanced Study during the writing of this paper.

\section{Local well-posedness}
The local well-posedness result of \cite{LS} may be proved via the Strichartz estimates of \cite{Stri}.
\begin{theorem}\label{t6.1}
Let $I \subset \mathbf{R}$, $t_{0} \in I$, be an interval and let $u$ solve the linear wave equation
\begin{equation}\label{6.1}
u_{tt} - \Delta u = F, \hspace{5mm} u(t_{0}) = u_{0}, \hspace{5mm} u_{t}(t_{0}) = u_{1}.
\end{equation}
Then we have the estimates
\begin{equation}\label{6.2}
\aligned
\| u \|_{L_{t}^{p} L_{x}^{q}(I \times \mathbf{R}^{3})} + \| u \|_{L_{t}^{\infty} \dot{H}^{s}(I \times \mathbf{R}^{3})} + \| u_{t} \|_{L_{t}^{\infty} \dot{H}^{s - 1}(I \times \mathbf{R}^{3})} \\
\lesssim_{p, q, s, \tilde{p}, \tilde{q}} \| u_{0} \|_{\dot{H}^{s}(\mathbf{R}^{3})} + \| u_{1} \|_{\dot{H}^{s - 1}(\mathbf{R}^{3})} + \| F \|_{L_{t}^{\tilde{p}'} L_{x}^{\tilde{q}'}(I \times \mathbf{R}^{3})},
\endaligned
\end{equation}
whenever $s \geq 0$, $2 \leq p, \tilde{p} \leq \infty$, $2 \leq q, \tilde{q} < \infty$, and
\begin{equation}\label{6.3}
\frac{1}{p} + \frac{1}{q} \leq \frac{1}{2}, \hspace{5mm} \frac{1}{\tilde{p}} + \frac{1}{\tilde{q}} \leq \frac{1}{2}.
\end{equation}
\end{theorem}
\emph{Proof:} Theorem $\ref{t6.1}$ was proved for $p = q = 4$ in \cite{Stri} and then in \cite{GV} for a general choice of $(p, q)$. $\Box$\vspace{5mm}

To prove local well-posedness of $(\ref{1.1})$, $(\ref{6.2})$ when $p = q = 4$ will suffice. Indeed, $(\ref{6.2})$ implies that for any $I$,
\begin{equation}\label{6.4}
\| u \|_{L_{t,x}^{4}(I \times \mathbf{R}^{3})} \lesssim \| S(t)(u_{0}, u_{1}) \|_{L_{t,x}^{4}(I \times \mathbf{R}^{3})} + \| u \|_{L_{t,x}^{4}(I \times \mathbf{R}^{3})}^{3}.
\end{equation}
If $\| S(t)(u_{0}, u_{1}) \|_{L_{t,x}^{4}(I \times \mathbf{R}^{3})} \leq \epsilon$, then $(\ref{6.4})$ implies that $(\ref{1.1})$ is locally well-posed on the interval $I$.

For $\| u_{0} \|_{\dot{H}^{1/2}} + \| u_{1} \|_{\dot{H}^{-1/2}}$ sufficiently small, $(\ref{6.2})$ and $(\ref{6.4})$ imply that $(\ref{1.1})$ is well-posed on $I = \mathbf{R}$. For generic $(u_{0}, u_{1}) \in \dot{H}^{1/2} \times \dot{H}^{-1/2}$, the dominated convergence theorem and $(\ref{6.2})$ imply that for any fixed $(u_{0}, u_{1}) \in \dot{H}^{1/2} \times \dot{H}^{-1/2}$,
\begin{equation}\label{6.5}
\lim_{T \searrow 0} \| S(t)(u_{0}, u_{1}) \|_{L_{t,x}^{4}([-T, T] \times \mathbf{R}^{3})} = 0,
\end{equation}
which implies local well-posedness on some open interval $I$, where $0 \in I$.\vspace{5mm}

Equation $(\ref{6.4})$ also implies that $(\ref{1.1})$ is locally well-posed on an interval $I$ on which an a priori bound $\| u \|_{L_{t,x}^{4}(I \times \mathbf{R}^{3})} < \infty$ is obtained. This may be seen by partitioning $I$ into finitely many pieces $I_{j}$ on which $\| u \|_{L_{t,x}^{4}(I_{j} \times \mathbf{R}^{3})}$ is small, and then iterating local well-posedness arguments on each interval. This argument also shows that scattering is equivalent to $\| u \|_{L_{t,x}^{4}(\mathbf{R} \times \mathbf{R}^{3})} < \infty$.\vspace{5mm}

Strichartz estimates also yield perturbative results.
\begin{lemma}[Perturbation lemma]\label{l6.2}
Let $I \subset \mathbf{R}$ be a time interval. Let $t_{0} \in I$, $(u_{0}, u_{1}) \in \dot{H}^{1/2} \times \dot{H}^{-1/2}$ and some constants $M$, $A$, $A' > 0$. Let $\tilde{u}$ solve the equation
\begin{equation}\label{6.6}
(\partial_{tt} - \Delta) \tilde{u} = F(\tilde{u}) = e,
\end{equation}
on $I \times \mathbf{R}^{3}$, and also suppose $\sup_{t \in I} \| (\tilde{u}(t), \partial_{t} \tilde{u}(t)) \|_{\dot{H}^{1/2} \times \dot{H}^{-1/2}} \leq A$, $\| \tilde{u} \|_{L_{t,x}^{4}(I \times \mathbf{R}^{3})} \leq M$,
\begin{equation}\label{6.7}
\| (u_{0} - \tilde{u}(t_{0}), u_{1} - \partial_{t} \tilde{u}(t_{0})) \|_{\dot{H}^{1/2} \times \dot{H}^{-1/2}} \leq A',
\end{equation}
and
\begin{equation}\label{6.8}
\| e \|_{L_{t,x}^{4/3}(I \times \mathbf{R}^{3})} + \| S(t - t_{0})(u_{0} - \tilde{u}(t_{0}), u_{1} - \partial_{t} \tilde{u}(t_{0})) \|_{L_{t,x}^{4}(I \times \mathbf{R}^{3})} \leq \epsilon.
\end{equation}
Then there exists $\epsilon_{0}(M, A, A')$ such that if $0 < \epsilon < \epsilon_{0}$ then there exists a solution to $(\ref{1.1})$ on $I$ with $(u(t_{0}), \partial_{t} u(t_{0})) = (u_{0}, u_{1})$, $\| u \|_{L_{t,x}^{4}(I \times \mathbf{R}^{3})} \leq C(M, A, A')$, and for all $t \in I$,
\begin{equation}\label{6.9}
\| (u(t), \partial_{t} u(t)) - (\tilde{u}(t), \partial_{t} \tilde{u}(t)) \|_{\dot{H}^{1/2} \times \dot{H}^{-1/2}} \leq C(A, A', M)(A' + \epsilon).
\end{equation}
\end{lemma}
\emph{Proof:} The method of proof is by now fairly well-known. See for example lemma $2.20$ of \cite{KM}. $\Box$\vspace{5mm}

The proof of Theorem $\ref{t1.2}$ also utilizes some additional Strichartz estimates that only appear for radially symmetric data. First, \cite{KlMa} proved that the endpoint case of Theorem $\ref{t6.1}$ also holds.
\begin{theorem}\label{t6.3}
For $(u_{0}, u_{1})$ radially symmetric, and $u$ solves $(\ref{6.1})$ with $F = 0$,
\begin{equation}\label{6.10}
\| u \|_{L_{t}^{2} L_{x}^{\infty}(\mathbf{R} \times \mathbf{R}^{3})} \lesssim \| u_{0} \|_{\dot{H}^{1}(\mathbf{R}^{3})} + \| u_{1} \|_{L^{2}(\mathbf{R}^{3})}.
\end{equation}
\end{theorem}

Additionally, the proof will rely very heavily on the estimates of \cite{Sterb} for radially symmetric initial data, extending the range of $(p, q)$ in $(\ref{6.3})$ for radial initial data.
\begin{theorem}\label{t6.4}
Let $(u_{0}, u_{1})$ be spherically symmetric, and suppose $u$ solves $(\ref{6.1})$ with $F = 0$. Then if $q > 4$ and
\begin{equation}\label{6.11}
\frac{1}{2} + \frac{3}{q} = \frac{3}{2} - s,
\end{equation}
then
\begin{equation}\label{6.12}
\| u \|_{L_{t}^{2} L_{x}^{q}(\mathbf{R} \times \mathbf{R}^{3})} \lesssim \| u_{0} \|_{\dot{H}^{s}(\mathbf{R}^{3})} + \| u_{1} \|_{\dot{H}^{s - 1}(\mathbf{R}^{3})}.
\end{equation}
\end{theorem}

\section{Virial identities for the wave equation}
The proof of Theorem $\ref{t1.2}$ will also use some weighted Strichartz-type estimates. These estimates could actually be proved using Proposition $3.5$ of \cite{Sterb} after making a Bessel function-type reduction from three dimensions to two dimensions using radial symmetry.

However, these estimates will instead be proved using virial identities. There are at least two reasons for doing this. The first is that, in the author's opinion, the exposition is cleaner and more readable using virial identities. The second reason is that many of the computations may be applied equally well to defocusing problems as to linear problems.\vspace{5mm}

Suppose $u$ solves the equation
\begin{equation}\label{2.0}
u_{tt} - \Delta + \mu u^{3} = 0, \qquad u(0,x) = u_{0}, \qquad u_{t}(0,x) = u_{1},
\end{equation}
where $\mu = 0, 1$. The case when $\mu = 0$ is a solution to the linear wave equation and $\mu = 1$ is the defocusing nonlinear wave equation $(\ref{1.1})$.
\begin{theorem}\label{t2.3}
If $u$ solves $(\ref{1.1})$ on an interval $[0, T]$, then
\begin{equation}\label{2.19}
\int_{0}^{T} \int \frac{\mu}{|x|} u^{4} dx dt \lesssim \| u \|_{L_{t}^{\infty} \dot{H}^{1}([0, T] \times \mathbf{R}^{3})} \| u_{t} \|_{L_{t}^{\infty} L_{x}^{2}([0, T] \times \mathbf{R}^{3})},
\end{equation}
\begin{equation}\label{2.20}
\sup_{R > 0} \frac{1}{R^{3}} \int_{0}^{T} \int_{|x| \leq R} u^{2} dx dt \lesssim \| u \|_{L_{t}^{\infty} \dot{H}^{1}([0, T] \times \mathbf{R}^{3})} \| u_{t} \|_{L_{t}^{\infty} L_{x}^{2}([0, T] \times \mathbf{R}^{3})},
\end{equation}
and
\begin{equation}\label{2.21}
\sup_{R > 0} \frac{1}{R} \int_{0}^{T} \int_{|x| \leq R} [|\nabla u|^{2} + u_{t}^{2}] dx dt \lesssim \| u \|_{L_{t}^{\infty} \dot{H}^{1}([0, T] \times \mathbf{R}^{3})} \| u_{t} \|_{L_{t}^{\infty} L_{x}^{2}([0, T] \times \mathbf{R}^{3})}.
\end{equation}
\end{theorem}
\noindent \emph{Proof:} Define the generic Morawetz potential
\begin{equation}\label{2.22}
M(t) = \int u_{t} a(|x|) x \cdot \nabla u + \int u_{t} a(|x|) u.
\end{equation}
Computing the time derivative,
\begin{equation}\label{2.23}
\aligned
\frac{d}{dt} M(t) = \int u_{t} a(|x|) x \cdot \nabla u_{t} + \int u_{t}^{2} a(|x|) \\ \int \Delta u a(|x|) x \cdot \nabla u + \int \Delta u a(|x|) u \\ - \mu \int u^{3} a(|x|) x \cdot \nabla u - \mu \int u^{3} a(|x|) u.
\endaligned
\end{equation}
Integrating by parts,
\begin{equation}\label{2.24}
\aligned
\frac{d}{dt} M(t) = -\frac{1}{2} \int [a(|x|) + a'(|x|) |x|] u_{t}^{2} - \frac{1}{2} \int [a(|x|) + a'(|x|) |x|] |\nabla u|^{2} \\
+ \int a'(|x|) |x| [|\nabla u|^{2} - |\partial_{r} u|^{2}] + \frac{1}{2} \int u^{2} \Delta a(|x|) \\ - \frac{\mu}{4} \int a(|x|) u^{4} + \frac{\mu}{4} \int a'(|x|) |x| u^{4}.
\endaligned
\end{equation}
If we choose $a(|x|) = \frac{1}{|x|}$, then
\begin{equation}\label{2.25}
a(|x|) + a'(|x|) |x| = 0.
\end{equation}
When $u$ is radial, $|\nabla u|^{2} - |\partial_{r} u|^{2} = 0$. For a general $u$,
\begin{equation}\label{2.26}
|\nabla u|^{2} - |\partial_{r} u|^{2} \geq 0,
\end{equation}
so since $a'(|x|) \leq 0$,
\begin{equation}\label{2.27}
a'(|x|) |x| [|\nabla u|^{2} - |\partial_{r} u|^{2}] \leq 0.
\end{equation}
Also, by direct calculation, $\Delta \frac{1}{|x|} = - 2 \pi \delta(x)$, so when $a(|x|) = \frac{1}{|x|}$,
\begin{equation}\label{2.28}
\frac{d}{dt} M(t) \leq -\pi u(t,0)^{2} - \frac{\mu}{2} \int \frac{1}{|x|} u^{4} dx.
\end{equation}
Now by Hardy's inequality, when $a(x) = \frac{1}{|x|}$,
\begin{equation}\label{2.29}
|M(t)| \lesssim \| u_{t} \|_{L^{2}} \| \nabla u \|_{L^{2}}.
\end{equation}
Therefore,
\begin{equation}\label{2.30}
\int_{0}^{T} u(t,0)^{2} dt + \int_{0}^{T} \int \frac{\mu}{|x|} u^{4} dx dt \lesssim \| u_{t} \|_{L_{t}^{\infty} L_{x}^{2}} \| \nabla u \|_{L_{t}^{\infty} L_{x}^{2}}.
\end{equation}
This takes care of $(\ref{2.19})$. 

Replacing $a(|x|)$ by $a(|x - y|)$ and $x$ with $x - y$, $(\ref{2.30})$ implies
\begin{equation}\label{2.31}
\aligned
\frac{1}{R^{3}} \int_{0}^{T} \int_{|y| \leq R} u(t,y)^{2} dy dt + \frac{1}{R^{3}} \int_{|y| \leq R} \int \frac{\mu}{|x - y|} u(t,x)^{4} dx dy \\ \lesssim \| u_{t} \|_{L_{t}^{\infty} L_{x}^{2}} \| \nabla u \|_{L_{t}^{\infty} L_{x}^{2}},
\endaligned
\end{equation}
which takes care of $(\ref{2.20})$.

To prove $(\ref{2.21})$, choose a smooth function $\chi : [0, \infty) \rightarrow [0, \infty)$ satisfying $\chi(|x|) = 1$ for $|x| \leq 1$, $\chi(|x|) = \frac{3}{2 |x|}$ for $|x| \geq 2$, and such that
\begin{equation}
\chi(|x|) + \chi'(|x|) |x| = \phi(|x|)
\end{equation}
is a smooth function, $\phi(|x|) \geq 0$, $\phi(|x|) = 1$ for $|x| \leq 1$, and $\phi(|x|)$ is supported on $|x| \leq 2$. Take $a(|x|) = \frac{1}{R} \chi(\frac{|x|}{R})$.
\begin{equation}\label{2.32}
a(|x|) + a'(|x|) |x| = \frac{1}{R} \chi(\frac{|x|}{R}) + \frac{1}{R} \chi'(\frac{|x|}{R}) \frac{|x|}{R} = \frac{1}{R} \phi(\frac{|x|}{R}).
\end{equation}
Therefore,
\begin{equation}\label{2.33}
\frac{d}{dt} M(t) = -\frac{1}{2R} \int \phi(\frac{|x|}{R}) [u_{t}^{2} + |\nabla u|^{2}] - \frac{\mu}{4R} \int a(\frac{|x|}{R}) u^{4} + \frac{1}{2R} \int u^{2} \Delta a(\frac{|x|}{R}).
\end{equation}
Now, since $a(|x|) = \frac{3}{2} \frac{1}{|x|}$ when $|x| \geq 2$, $\Delta a(|x|)$ is supported on $|x| \leq 2$.  Therefore,
\begin{equation}\label{2.34}
\frac{1}{2R} \int u^{2} \Delta a(\frac{|x|}{R}) \lesssim \sup_{R > 0} \frac{1}{R^{3}} \int_{|x| \leq R} u^{2}.
\end{equation} 
Also, $a(|x|) \lesssim \frac{1}{|x|}$, so again by Hardy's inequality,
\begin{equation}\label{2.33.1}
|M(t)| \lesssim \| v_{t} \|_{L^{2}} \| \nabla v \|_{L^{2}}.
\end{equation}
Plugging $(\ref{2.31})$ and $(\ref{2.34})$ into $(\ref{2.33})$ proves $(\ref{2.21})$. $\Box$

\begin{corollary}\label{c2.4}
If $u$ is an approximate solution to the cubic wave equation,
\begin{equation}\label{2.35}
u_{tt} - \Delta u + u^{3} = F,
\end{equation}
then
\begin{equation}\label{2.36}
\aligned
\frac{d}{dt} [\int u_{t} \frac{x}{|x|} \cdot \nabla u + \int u_{t} \frac{1}{|x|} u] \leq -2 \pi u(t,0)^{2} - \frac{1}{2} \int \frac{1}{|x|} u^{4} \\
+ \int F \frac{x}{|x|} \cdot \nabla u + \int F \frac{1}{|x|} u,
\endaligned
\end{equation}
\begin{equation}\label{2.37}
\aligned
\frac{d}{dt} \frac{1}{R^{3}} [\int_{|y| \leq R} \int u_{t} \frac{x}{|x|} \cdot \nabla u + \int_{|y| \leq R} \int u_{t} \frac{1}{|x|} u] \\ \leq -\pi \frac{1}{R^{3}} \int_{|y| \leq R} u(t,y)^{2} - \frac{1}{2} \frac{1}{R^{3}} \int_{|y| \leq R} \int \frac{1}{|x - y|} u^{4} \\
+ \frac{1}{R^{3}} \int_{|y| \leq R} \int F \frac{x - y}{|x - y|} \cdot \nabla u + \frac{1}{R^{3}} \int_{|y| \leq R} \int F \frac{1}{|x - y|} u,
\endaligned
\end{equation}
and
\begin{equation}\label{2.38}
\aligned
\frac{d}{dt} [\frac{1}{R} \int u_{t} \chi(\frac{|x|}{R}) x \cdot \nabla u + \frac{1}{R} \int u_{t} \chi(\frac{|x|}{R}) u] \leq -\frac{1}{2R} \int \phi(\frac{|x|}{R}) [u_{t}^{2} + |\nabla u|^{2}] \\
- \frac{1}{4R} \int \chi(\frac{|x|}{R}) u^{4} + \frac{1}{4R} \int \chi'(\frac{|x|}{R}) \frac{|x|}{R} u^{4} \\ + \frac{1}{2R} \int u^{2} \Delta a(\frac{|x|}{R}) + \frac{1}{R} \int F \chi(\frac{|x|}{R}) x \cdot \nabla u + \frac{1}{R} \int F \chi(\frac{|x|}{R}) u.
\endaligned
\end{equation}
\end{corollary}

Theorem $\ref{t2.3}$ also gives some nice estimates for the linear wave equation $(\mu = 0)$.
\begin{corollary}\label{c2.1}
For any $j \in \mathbf{Z}$, let $w$ be the solution to the linear wave equation
\begin{equation}
\partial_{tt} w - \Delta w = 0, \qquad w(0,x) = P_{j} u_{0}, \qquad w_{t}(0,x) = P_{j} u_{1}.
\end{equation}
Then for any $2 < p < \infty$,
\begin{equation}\label{2.1}
\| |x|^{1/2} w \|_{L_{t}^{p} L_{x}^{\infty}(\mathbf{R} \times \mathbf{R}^{3})} \lesssim \| P_{j} u_{0} \|_{\dot{H}^{1/p'}(\mathbf{R}^{3})} + \| P_{j} u_{1} \|_{\dot{H}^{1/p' - 1}(\mathbf{R}^{3})},
\end{equation}
where $\frac{1}{p'} = 1 - \frac{1}{p}$ is the Lebesgue dual of $p$. Also, for $p = 2$, for any $0 < R < 1$, $1 < R_{1} < \infty$,
\begin{equation}\label{2.2}
\| |x|^{1/2} w \|_{L_{t,x}^{2}(\mathbf{R} \times \{ x : R \leq |x| \leq R_{1} \})}^{2} \lesssim (\ln(R_{1}) - \ln(R) + 1) [\| P_{j} u_{0} \|_{\dot{H}^{1/2}(\mathbf{R}^{3})}^{2} + \| P_{j} u_{1} \|_{\dot{H}^{-1/2}(\mathbf{R}^{3})}^{2}].
\end{equation}
\end{corollary}
\emph{Proof:} Let $\psi$ be a smooth, radial function supported on an annulus, $\psi(r) = 1$ for $1 \leq r \leq 2$, and $\psi(r)$ is supported on $\frac{1}{2} \leq r \leq 4$. By Bernstein's inequality,
\begin{equation}\label{2.11}
\| P_{k}(\psi(\frac{r}{R}) w) \|_{L^{2}} \lesssim 2^{-k} \| \partial_{r}(\psi(\frac{r}{R}) w) \|_{L^{2}} + 2^{-k} R^{-1} \| \psi'(\frac{r}{R}) w \|_{L^{2}}.
\end{equation}
Therefore, by $(\ref{2.20})$, $(\ref{2.21})$, and the radial Sobolev embedding theorem,
\begin{equation}\label{2.12}
\sum_{k \geq j - 3} \| P_{k} (\psi(\frac{r}{R}) w) \|_{L_{t}^{2} L_{x}^{\infty}} \lesssim 2^{-j/2} R^{-1/2} (\| P_{j} u_{0} \|_{\dot{H}^{1}} + \| P_{j} u_{1} \|_{L^{2}}).
\end{equation}
Next, by the Fourier support properties of $w$,
\begin{equation}\label{2.13}
\| P_{\leq j - 3} (\psi(\frac{r}{R}) w) \|_{L^{\infty}} \lesssim 2^{-j} R^{-1} \| w \|_{L^{\infty}}.
\end{equation}
Combining $(\ref{2.13})$ with $(\ref{6.10})$,
\begin{equation}\label{2.14}
\| P_{\leq j - 3} (\psi(\frac{r}{R}) w) \|_{L_{t}^{2} L_{x}^{\infty}} \lesssim 2^{-j} R^{-1}  (\| P_{j} u_{0} \|_{\dot{H}^{1}} + \| P_{j} u_{1} \|_{L^{2}}).
\end{equation}
Then when $R \geq 2^{-j}$,
\begin{equation}
\| P_{\leq j - 3} (\psi(\frac{r}{R}) w) \|_{L_{t}^{2} L_{x}^{\infty}} \lesssim 2^{-j/2} R^{-1/2}  (\| P_{j} u_{0} \|_{\dot{H}^{1}} + \| P_{j} u_{1} \|_{L^{2}}).
\end{equation}
Finally, when $R \leq 2^{-j}$, a straightforward application of the endpoint Strichartz estimate yields
\begin{equation}\label{2.15}
\| \psi(\frac{r}{R}) w \|_{L_{t}^{2} L_{x}^{\infty}} \lesssim (\| P_{j} u_{0} \|_{\dot{H}^{1}} + \| P_{j} u_{1} \|_{L^{2}}) \lesssim R^{-1/2} 2^{-j/2} (\| P_{j} u_{0} \|_{\dot{H}^{1}} + \| P_{j} u_{1} \|_{L^{2}}).
\end{equation}
Since there are $\lesssim \ln(R_{1}) - \ln(R) + 1$ dyadic annuli overlapping $R \leq |x| \leq R_{1}$, $(\ref{2.12})$--$(\ref{2.15})$ directly yields $(\ref{2.2})$.

To prove $(\ref{2.1})$, interpolating $(\ref{2.15})$ with the radial Sobolev embedding theorem, for any $2 < p < \infty$,
\begin{equation}\label{2.16}
\aligned
\| \psi(\frac{r}{R}) w \|_{L_{t}^{p} L_{x}^{\infty}} \lesssim \| \psi(\frac{r}{R}) w \|_{L_{t}^{2} L_{x}^{\infty}}^{2/p} \| \psi(\frac{r}{R}) w \|_{L_{t,x}^{\infty}}^{1 - 2/p} \\ \lesssim R^{-1/2} R^{-\frac{1}{2}(1 - \frac{2}{p})} (\| P_{j} u_{0} \|_{\dot{H}^{1/2}} + \| P_{j} u_{1} \|_{\dot{H}^{-1/2}}),
\endaligned
\end{equation}
which directly implies
\begin{equation}\label{2.17}
\| |x|^{1/2} w \|_{L_{t}^{p} L_{x}^{\infty}(\mathbf{R} \times \{ x : |x| \geq 2^{-j} \})} \lesssim (\| P_{j} u_{0} \|_{\dot{H}^{1/p'}} + \| P_{j} u_{1} \|_{\dot{H}^{1/p' - 1}}).
\end{equation}
Meanwhile, by $(\ref{6.10})$ and the Sobolev embedding theorem,
\begin{equation}\label{2.18}
\| |x|^{1/2} w \|_{L_{t}^{p} L_{x}^{\infty}(\mathbf{R} \times \{ |x| \leq 2^{-j} \})} \lesssim 2^{-j/2} \| w \|_{L_{t}^{2} L_{x}^{\infty}}^{2/p} \| w \|_{L_{t,x}^{\infty}}^{1 - 2/p} \lesssim (\| P_{j} u_{0} \|_{\dot{H}^{1/p}} + \| P_{j} u_{1} \|_{\dot{H}^{1/p - 1}}).
\end{equation}
This finally proves the theorem. $\Box$\medskip

\noindent \textbf{Remark:} Also observe that by the radial Sobolev embedding theorem, Corollary $\ref{c2.1}$ implies
\begin{equation}\label{2.18.1}
\| w \|_{L_{t}^{2} L_{x}^{\infty}([0, T] \times \{ |x| \geq R \})}^{2} \lesssim (1 + \ln(T) - \ln(R)) [\| P_{j} u_{0} \|_{\dot{H}^{1/2}} + \| P_{j} u_{1} \|_{\dot{H}^{-1/2}}].
\end{equation}

The virial identities in Theorem $\ref{t2.3}$ commute very well with Littlewood--Paley projections.
\begin{lemma}\label{l3.1}
For any $j$,
\begin{equation}\label{3.18}
\int \frac{1}{|x|} |P_{\leq j} v|^{4} dx + \int \frac{1}{|x|} |P_{\geq j} v|^{4} dx \lesssim \int \frac{1}{|x|} |v|^{4} dx.
\end{equation}
\end{lemma}
\noindent \emph{Proof:} Let $\psi$ be the Littlewood--Paley kernel.
\begin{equation}\label{2.19.1}
\frac{1}{|x|^{1/4}} P_{\leq j} v(x) = \frac{1}{|x|^{1/4}} \int 2^{3j} \psi(2^{j}(x - y)) v(y) dy.
\end{equation}
When $|y| \lesssim |x|$,
\begin{equation}\label{2.20.1}
\frac{1}{|x|^{1/4}} 2^{3j} \psi(2^{j}(x - y)) \lesssim 2^{3j} \psi(2^{j}(x - y)) \frac{1}{|y|^{1/4}}.
\end{equation}
When $|y| \gg |x|$ and $|x| \geq 2^{-j}$, since $\psi$ is rapidly decreasing, for any $N$,
\begin{equation}\label{2.21.1}
\aligned
\frac{1}{|x|^{1/4}} 2^{3j} \psi(2^{j}(x - y)) \lesssim_{N} \frac{1}{|x|^{1/4}} \frac{2^{3j}}{(1 + 2^{j} |x - y|)^{N}} \\ \lesssim \frac{1}{|x|^{1/4} 2^{j} |y|} \frac{2^{3j}}{(1 + 2^{j} |x - y|)^{N - 1}} \lesssim \frac{1}{|y|^{1/4}} \frac{2^{3j}}{(1 + 2^{j}|x - y|)^{N - 1}}.
\endaligned
\end{equation}
Combining $(\ref{2.20.1})$ and $(\ref{2.21.1})$,
\begin{equation}\label{2.22}
\| \frac{1}{|x|^{1/4}} |P_{\leq j} v| \|_{L^{4}(|x| \geq 2^{-j})} \lesssim \| \frac{1}{|x|^{1/4}} v \|_{L^{4}(\mathbf{R}^{3})}.
\end{equation}
When $|y| \gg |x|$ and $|x| \leq 2^{-j}$, since $\psi$ is rapidly decreasing, for any $N$,
\begin{equation}\label{2.23}
\aligned
\frac{1}{|x|^{1/4}} 2^{3j} \psi(2^{j}(x - y)) \lesssim_{N} \frac{1}{|x|^{1/4}} \frac{2^{3j}}{(1 + 2^{j} |x - y|)^{N}} \\ \lesssim \frac{1}{|x|^{1/4}} \frac{2^{3j}}{(1 + 2^{j} |x - y|)^{N - 1/4}} \frac{1}{2^{j/4} |y|^{1/4}}.
\endaligned
\end{equation}
\begin{equation}\label{2.24}
\| \frac{2^{11j/4}}{(1 + 2^{j}|x - y|)^{N}} \|_{L^{4/3}(\mathbf{R}^{3})} \lesssim 2^{j/2},
\end{equation}
so by $(\ref{2.20.1})$, $(\ref{2.24})$, Young's inequality, and H{\"o}lder's inequality,
\begin{equation}\label{2.25}
\| \frac{1}{|x|^{1/4}} |P_{\leq j} v| \|_{L^{4}(|x| \leq 2^{-j})} \lesssim \| \frac{1}{|x|^{1/4}} v \|_{L^{4}(\mathbf{R}^{3})}.
\end{equation}
This proves $(\ref{3.18})$. $\Box$
\begin{lemma}\label{l2.2}
\begin{equation}\label{2.26}
\aligned
\| P_{\geq j} v \|_{L^{4}(|x| \leq \frac{R}{2})}^{2} \lesssim \| P_{\geq j} v \|_{L^{3}} [\| \nabla v \|_{L^{2}(|x| \leq R)} + \frac{1}{R} \| v \|_{L^{2}(|x| \leq R)}] \\ + 2^{-j/2} (\int \frac{1}{|x|} v^{4})^{1/2}.
\endaligned
\end{equation}
\end{lemma}
\emph{Proof:} Let $\phi \in C_{0}^{\infty}(\mathbf{R}^{3})$ be supported on $|x| \leq 1$ and $\phi(x) = 1$ for $|x| \leq \frac{1}{2}$. By H{\"o}lder's inequality,
\begin{equation}\label{2.27}
\| P_{\geq j} v \|_{L^{4}(|x| \leq \frac{R}{2})}^{2} \leq \| \phi(\frac{x}{R}) (P_{\geq j} v) \|_{L^{4}(\mathbf{R}^{3})}^{2}.
\end{equation}
Then, by the triangle inequality, H{\"o}lder's inequality, and the Cauchy--Schwartz inequality,
\begin{equation}\label{2.28}
\aligned
\| \phi(\frac{x}{R}) (P_{\geq j} v) \|_{L^{4}(\mathbf{R}^{3})}^{2} \leq \| \phi(\frac{x}{R}) (P_{\geq j} v) \cdot P_{\geq j} (\phi(\frac{x}{R}) v) \|_{L^{2}(\mathbf{R}^{3})} \\ + \| \phi(\frac{x}{R}) (P_{\geq j} v) \cdot [\phi(\frac{x}{R}), P_{\geq j}] v \|_{L^{2}(\mathbf{R}^{3})} \leq \| P_{\geq j} \phi(\frac{x}{R}) v \|_{L^{6}(\mathbf{R}^{3})} \| P_{\geq j} v \|_{L^{3}(\mathbf{R}^{3})} \\ + \frac{1}{2} \| \phi(\frac{x}{R}) (P_{\geq j} v) \|_{L^{4}(\mathbf{R}^{3})}^{2} + \frac{1}{2} \| [\phi(\frac{x}{R}), P_{\geq j}] v \|_{L^{4}(\mathbf{R}^{3})}^{2},
\endaligned
\end{equation}
where
\begin{equation}\label{2.29.1}
[\phi(\frac{x}{R}), P_{\geq j}] v = \phi(\frac{x}{R}) (P_{\geq j} v) - P_{\geq j} (\phi(\frac{x}{R}) v).
\end{equation}
By the Littlewood--Paley theorem,
\begin{equation}\label{2.30.1}
\| \phi(\frac{x}{R}) (P_{\geq j} v) \|_{L^{4}(\mathbf{R}^{3})}^{2} \lesssim \| \phi(\frac{x}{R}) v \|_{L^{6}(\mathbf{R}^{3})} \| P_{\geq j} v \|_{L^{3}(\mathbf{R}^{3})} + \| [\phi(\frac{x}{R}), P_{\geq j}] v \|_{L^{4}(\mathbf{R}^{3})}^{2},
\end{equation}
and by the Sobolev embedding theorem,
\begin{equation}\label{2.31.1}
\| \phi(\frac{x}{R}) v \|_{L^{6}(\mathbf{R}^{3})} \lesssim \| \nabla (\phi(\frac{x}{R}) v) \|_{L^{2}(\mathbf{R}^{3})} \lesssim \frac{1}{R} \| v \|_{L^{2}(|x| \leq R)} + \| \nabla v \|_{L^{2}(|x| \leq R)}.
\end{equation}
This takes care of the first term on the right hand side of $(\ref{2.30.1})$.

To handle the commutator, observe that
\begin{equation}\label{2.32}
[\phi(\frac{x}{R}), P_{\geq j}] = -[P_{\leq j}, \phi(\frac{x}{R})].
\end{equation}
Then compute
\begin{equation}\label{2.33.1.1}
[P_{\leq j}, \phi(\frac{x}{R})] v = 2^{3j} \int \psi(2^{j}(x - y)) [\phi(\frac{y}{R}) - \phi(\frac{x}{R})] v(y) dy.
\end{equation}
When $|y| \gg |x|$, the kernel
\begin{equation}\label{2.34.1}
\aligned
2^{3j} \psi(2^{j}(x - y)) [\phi(\frac{y}{R}) - \phi(\frac{x}{R})] \lesssim_{N} \frac{2^{3j}}{(1 + 2^{j} |x - y|)^{N}} \\ \lesssim 2^{-j/4} \frac{2^{3j}}{(1 + 2^{j} |x - y|)^{N - 1/4}} \frac{1}{|y|^{1/4}}.
\endaligned
\end{equation}
When $|y| \lesssim |x|$ and $|x| \leq R$, by the fundamental theorem of calculus,
\begin{equation}\label{2.35}
\aligned
2^{3j} \psi(2^{j}(x - y)) [\phi(\frac{y}{R}) - \phi(\frac{x}{R})] \lesssim_{N} \frac{2^{3j}}{(1 + 2^{j} |x - y|)^{N}} \frac{|x - y|^{1/4}}{R^{1/4}} \\ \lesssim 2^{-j/4} \frac{2^{3j}}{(1 + 2^{j} |x - y|)^{N - 1/4}} \cdot \frac{1}{|y|^{1/4}}.
\endaligned
\end{equation}
When $|y| \lesssim |x|$ and $|x| > R$, interpolating
\begin{equation}\label{2.36.1}
\aligned
2^{3j} \psi(2^{j}(x - y)) [\phi(\frac{y}{R}) - \phi(\frac{x}{R})] = 2^{3j} \psi(2^{j}(x - y)) \phi(\frac{y}{R}) \lesssim_{N} \frac{2^{3j}}{(1 + 2^{j}|x - y|)^{N}} \frac{R^{1/2}}{|y|^{1/2}}
\endaligned
\end{equation}
with the fact that
\begin{equation}\label{2.37.1}
\aligned
2^{3j} \psi(2^{j}(x - y)) [\phi(\frac{y}{R}) - \phi(\frac{x}{R})] \lesssim_{N} \frac{2^{3j}}{(1 + 2^{j} |x - y|)^{N}} \frac{|x - y|^{1/2}}{R^{1/2}} \\ \lesssim 2^{-j/4} \frac{2^{3j}}{(1 + 2^{j} |x - y|)^{N - 1/2}} \cdot \frac{1}{2^{j/2} R^{1/2}},
\endaligned
\end{equation}
implies
\begin{equation}\label{2.38.1}
2^{3j} \psi(2^{j}(x - y)) [\phi(\frac{y}{R}) - \phi(\frac{x}{R})] \lesssim_{N} 2^{-j/4} \frac{2^{3j}}{(1 + 2^{j}|x - y|)^{N}} \frac{1}{|y|^{1/4}}.
\end{equation}
The kernel estimates $(\ref{2.34})$, $(\ref{2.35})$, and $(\ref{2.38})$ imply that
\begin{equation}\label{2.39}
\| [\phi(\frac{x}{R}), P_{\geq j}] v \|_{L^{4}(\mathbf{R}^{3})} \lesssim 2^{-j/4} \| \frac{1}{|x|^{1/4}} v \|_{L^{4}(\mathbf{R}^{3})},
\end{equation}
proving Lemma $\ref{l2.2}$. $\Box$

\section{Global well-posedness}
To prove global well-posedness of $(\ref{1.1})$ using the Fourier truncation method, decompose the initial data into a finite energy piece and a small data piece, $u_{0} = v_{0} + w_{0}$ and $u_{1} = v_{1} + w_{1}$, where
\begin{equation}
E(v_{0}, v_{1}) = \frac{1}{2} \int |\nabla v_{0}|^{2} dx + \frac{1}{2} \int |v_{1}|^{2} dx + \frac{1}{4} \int |v_{0}|^{4} dx < \infty,
\end{equation}
and
\begin{equation}
\| w_{0} \|_{\dot{H}^{1/2}} + \| w_{1} \|_{\dot{H}^{-1/2}} \ll 1.
\end{equation}
A solution $u$ to $(\ref{1.1})$ may then be decomposed into $u = w + v$, where $w$ solves
\begin{equation}\label{3.2}
w_{tt} - \Delta w + w^{3} = 0, \qquad w(0,x) = w_{0}, \qquad w_{t}(0,x) = w_{1},
\end{equation}
and $v$ solves
\begin{equation}\label{3.3}
v_{tt} - \Delta v + v^{3} + 3 v^{2} w + 3 v w^{2} = 0, \qquad v(0,x) = v_{0}, \qquad v_{t}(0,x) = v_{1}.
\end{equation}
If $\| w_{0} \|_{\dot{H}^{1/2}} + \| w_{1} \|_{\dot{H}^{-1/2}} < \epsilon$ for some $\epsilon > 0$ sufficiently small, then the small data arguments in $(\ref{6.4})$ implies that $(\ref{3.2})$ is globally well-posed, and moreover, by Theorem $\ref{t6.4}$,
\begin{equation}\label{3.4.1}
\aligned
\| w \|_{L_{t}^{2} L_{x}^{6}(\mathbf{R} \times \mathbf{R}^{3})} + \| |\nabla|^{1/10} w \|_{L_{t}^{2} L_{x}^{5}(\mathbf{R} \times \mathbf{R}^{3})} \\ + \| |\nabla|^{1/6} w \|_{L_{t}^{6} L_{x}^{3}(\mathbf{R} \times \mathbf{R}^{3})} + \| w \|_{L_{t,x}^{4}(\mathbf{R} \times \mathbf{R}^{3})} \lesssim \epsilon, \\ \| w^{3} \|_{L_{t}^{1} L_{x}^{3/2}(\mathbf{R} \times \mathbf{R}^{3})} \lesssim \epsilon^{3}.
\endaligned
\end{equation}
Following $(\ref{1.4})$, let $E(t)$ denote the energy of $v$, where
\begin{equation}\label{3.4}
E(t) = \frac{1}{2} \int v_{t}(t,x)^{2} dx + \frac{1}{2} \int |\nabla v(t,x)|^{2} dx + \frac{1}{4} \int v(t,x)^{4} dx.
\end{equation}

To prove global well-posedness it suffices to prove that $E(t) < \infty$ for all $t \in \mathbf{R}$.
\begin{theorem}\label{t3.1}
$(\ref{1.1})$ is locally well-posed on the time interval $[-\frac{c}{E(0)}, \frac{c}{E(0)}]$ for some fixed $c > 0$ sufficiently small.
\end{theorem}
\emph{Proof:} To simplify notation let $I = [-\frac{c}{E(0)}, \frac{c}{E(0)}]$. By Theorem $\ref{t1.1}$, $(\ref{1.1})$ has a solution for initial data $(v_{0}, v_{1})$, and moreover, by conservation of energy,
\begin{equation}\label{3.6}
\| v \|_{L_{t,x}^{4}(I \times \mathbf{R}^{3})}^{4} \lesssim |I| E(0) \leq c.
\end{equation}
Therefore, for $c > 0$ sufficiently small, independent of $E(0)$,
\begin{equation}\label{3.7}
\tilde{u}_{tt} - \Delta \tilde{u} + v^{3} + w^{3} = 0, \qquad \tilde{u}(0,x) = u_{0}, \qquad \tilde{u}_{t}(0,x) = u_{1},
\end{equation}
has a solution satisfying $\| \tilde{u} \|_{L_{t,x}^{4}(I \times \mathbf{R}^{3})} \ll 1$. Applying the perturbation lemma (Lemma $\ref{l6.2})$ completes the proof of Theorem $\ref{t3.1}$. $\Box$
\begin{theorem}\label{t3.2}
Equation $(\ref{1.1})$ is globally well-posed for radial $(u_{0}, u_{1}) \in \dot{H}^{1/2}(\mathbf{R}^{3}) \times \dot{H}^{-1/2}(\mathbf{R}^{3})$.
\end{theorem}
\noindent \emph{Proof:} To compute the time derivative of $E(t)$, by H{\"o}lder's inequality,
\begin{equation}\label{3.8}
\frac{d}{dt} E(t) = -3 \int v_{t} v^{2} w - 3 \int v_{t} v w^{2} \lesssim \| v_{t} \|_{L^{2}} \| v \|_{L^{6}}^{2} \| w \|_{L^{6}} + \| v_{t} \|_{L^{2}} \| v \|_{L^{6}} \| w \|_{L^{6}}^{2}.
\end{equation}
Therefore, by the Cauchy--Schwartz inequality,
\begin{equation}\label{3.10}
|\frac{d}{dt} E(t)| \lesssim E(t)^{2} + \| w \|_{L^{6}}^{2} E(t).
\end{equation}
If only the second term on the right hand side of $(\ref{3.10})$ were present, global boundedness of $E(t)$ would be an easy consequence of $(\ref{3.4.1})$ and Gronwall's inequality. However, the bound $|\frac{d}{dt} E(t)| \lesssim E(t)^{2}$ is not enough to exclude blow up in finite time. Instead, we will use a modification of $E(t)$, $\mathcal E(t)$, which has much better global derivative bounds, and satisfies $\mathcal E(t) \sim E(t)$.\vspace{5mm}

To simplify notation, rescale by $(\ref{1.2})$ so that
\begin{equation}\label{3.1}
\| P_{\geq 1} u_{0} \|_{\dot{H}^{1/2}(\mathbf{R}^{3})} + \| P_{\geq 1} u_{1} \|_{\dot{H}^{-1/2}(\mathbf{R}^{3})} < \epsilon,
\end{equation}
and then let $v_{0} = P_{\leq 1} u_{0}$ and $v_{1} = P_{\leq 1} u_{1}$. 

Following $(\ref{2.36})$, $(\ref{2.37})$, and $(\ref{2.38})$, let
\begin{equation}\label{3.1.1}
\aligned
M_{1}(t) = c_{1} \int v_{t} \frac{x}{|x|} \cdot \nabla v + c_{1} \int v_{t} \frac{1}{|x|} v, \\
M_{2}(t) = \frac{c_{2}}{R^{3}} \int_{|y| \leq 2R} \int v_{t} \frac{(x - y)}{|x - y|} \cdot \nabla v dx dy \\ + c_{2} \int_{|y| \leq 2R} \int v_{t} \frac{1}{|x - y|} v dx dy, \\
M_{3}(t) = \frac{c_{3}}{R}  \int v_{t} \chi(\frac{|x|}{R}) x \cdot \nabla v + \frac{c_{3}}{R} \int v_{t} \chi(\frac{|x|}{R}) v,
\endaligned
\end{equation}
where $c_{1}, c_{2}, c_{3} > 0$ are small constants and let
\begin{equation}
\mathcal E(t) = E(t) + M_{1}(t) + M_{2}(t) + M_{3}(t) + \int v^{3} w dx.
\end{equation}
Then by $(\ref{2.29})$, $(\ref{2.31})$, and $(\ref{2.33.1})$, and the Sobolev embedding theorem, which implies
\begin{equation}
\int v^{3} w dx \lesssim \| v \|_{L^{6}} \| w \|_{L^{3}} \| v \|_{L^{4}}^{2} \lesssim \epsilon E(t),
\end{equation}
we have
\begin{equation}\label{3.1.2}
\mathcal E(t) \sim E(t).
\end{equation}
Then by $(\ref{2.36})$, $(\ref{2.37})$, $(\ref{2.38})$, and $(\ref{3.8})$,
\begin{equation}\label{3.1.3}
\aligned
\frac{d}{dt} \mathcal E(t) \leq -c_{1} \pi v(t,0)^{2} - \frac{c_{2} \pi}{8 R^{3}} \int_{|y| \leq 2R} v(t,y)^{2} \\ - \frac{c_{1}}{2} \int \frac{1}{|x|} v^{4}  - \frac{c_{2}}{8 R^{3}}  \int_{|y| \leq 2R} \int \frac{1}{|x - y|} v^{4} \\ -\frac{c_{3}}{2R} \int \phi(\frac{|x|}{R}) [v_{t}^{2} + |\nabla v|^{2}] + \frac{c_{3}}{2R} \int v^{2} \Delta \chi(\frac{|x|}{R}) \\ - \frac{c_{3}}{4R} \int \chi(\frac{|x|}{R}) v^{4} + \frac{c_{3}}{4R} \int \chi'(\frac{|x|}{R}) \frac{|x|}{R} v^{4}  \\
+ \frac{d}{dt} \int v^{3} w dx + \int F v_{t} + c_{1} \int F \frac{x}{|x|} \cdot \nabla v + c_{1} \int F \frac{1}{|x|} v \\ + \frac{c_{2}}{8R^{3}} \int_{|y| \leq 2R} \int F \frac{(x - y)}{|x - y|} \cdot \nabla v + \frac{c_{2}}{8R^{3}} \int_{|y| \leq 2R} \int F \frac{1}{|x - y|} v \\ + \frac{c_{3}}{R} \int F \chi(\frac{|x|}{R}) x \cdot \nabla v + \frac{c_{3}}{R} \int F \chi(\frac{|x|}{R}) v,
\endaligned
\end{equation}
where $F = -3 v^{2} w - 3 v w^{2}$.

By the support properties of $\Delta \chi(\frac{|x|}{R})$, it is possible to choose $c_{2}, c_{3} > 0$ such that
\begin{equation}\label{3.1.4}
- \frac{c_{2} \pi}{8 R^{3}} \int_{|y| \leq 2R} v(t,y)^{2} + \frac{c_{3}}{2R} \int v^{2} \Delta \chi(\frac{|x|}{R}) \leq -\frac{c_{2}}{16 R^{3}} \int_{|y| \leq 2R} v(t,y)^{2}.
\end{equation}
Also, since $\chi'(\frac{|x|}{R}) \leq 0$,
\begin{equation}\label{3.1.5}
 \frac{c_{3}}{4R} \int \chi'(\frac{|x|}{R}) \frac{|x|}{R} v^{4} \leq 0.
\end{equation}
Therefore,
\begin{equation}\label{3.1.6}
\aligned
\frac{d}{dt} \mathcal E(t) + c_{1} \pi v(t,0)^{2} + \frac{c_{2} \pi}{16 R^{3}} \int_{|y| \leq 2R} v(t,y)^{2} \\ + \frac{c_{1}}{2} \int \frac{1}{|x|} v^{4}  + \frac{c_{2}}{8 R^{3}}  \int_{|y| \leq 2R} \int \frac{1}{|x - y|} v^{4} \\ +\frac{c_{3}}{2R} \int \phi(\frac{|x|}{R}) [v_{t}^{2} + |\nabla v|^{2}]  + \frac{c_{3}}{4R} \int \chi(\frac{|x|}{R}) v^{4} \\
\leq \frac{d}{dt} \int v^{3} w + \int F v_{t} + c_{1} \int F \frac{x}{|x|} \cdot \nabla v + c_{1} \int F \frac{1}{|x|} v \\ + \frac{c_{2}}{8R^{3}} \int_{|y| \leq 2R} \int F \frac{(x - y)}{|x - y|} \cdot \nabla v + \frac{c_{2}}{8R^{3}} \int_{|y| \leq 2R} \int F \frac{1}{|x - y|} v \\ + \frac{c_{3}}{R} \int F \chi(\frac{|x|}{R}) x \cdot \nabla v + \frac{c_{3}}{R} \int F \chi(\frac{|x|}{R}) v.
\endaligned
\end{equation}

By Hardy's inequality and the Sobolev embedding theorem,
\begin{equation}\label{3.26}
\int v^{2} w \frac{1}{|x|} v dx \lesssim (\int \frac{1}{|x|} v^{4} dx)^{1/2} \| \frac{1}{|x|^{1/2}} v \|_{L^{3}} \| w(t) \|_{L^{6}} \lesssim \delta (\int \frac{1}{|x|} v^{4} dx) + \frac{1}{\delta} E(t) \| w(t) \|_{L^{6}}^{2}.
\end{equation}
Also by H{\"o}lder's inequality and Hardy's inequality,
\begin{equation}\label{3.26.1}
\int v w^{2} \frac{1}{|x|} v \lesssim \| w \|_{L^{6}}^{2} \| \nabla v \|_{L^{2}} \| v \|_{L^{6}} \lesssim E(t) \| w \|_{L^{6}}^{2}.
\end{equation}
Therefore,
\begin{equation}\label{3.26.2}
\int F \frac{1}{|x|} v dx \lesssim \delta (\int \frac{1}{|x|} v^{4} dx) + \frac{1}{\delta} E(t) \| w(t) \|_{L^{6}}^{2}.
\end{equation}
Because $\chi(|x|) \lesssim \frac{1}{|x|}$, the same argument also implies
\begin{equation}\label{3.26.3}
\frac{1}{8 R^{3}} \int_{|y| \leq 2R} \int F \frac{1}{|x - y|} v + \frac{1}{R} \int F \chi(\frac{|x|}{R}) v \lesssim \delta (\int \frac{1}{|x|} v^{4} dx) + \frac{1}{\delta} E(t) \| w(t) \|_{L^{6}}^{2}.
\end{equation}
Therefore,
\begin{equation}\label{3.26.4}
\aligned
\frac{d}{dt} \mathcal E(t) + c_{1} \pi v(t,0)^{2} + \frac{c_{2} \pi}{16 R^{3}} \int_{|y| \leq 2R} v(t,y)^{2} \\ + \frac{c_{1}}{2} \int \frac{1}{|x|} v^{4}  + \frac{c_{2}}{8 R^{3}}  \int_{|y| \leq 2R} \int \frac{1}{|x - y|} v^{4} \\ +\frac{c_{3}}{2R} \int \phi(\frac{|x|}{R}) [v_{t}^{2} + |\nabla v|^{2}]  + \frac{c_{3}}{4R} \int \chi(\frac{|x|}{R}) v^{4} \\
- \frac{d}{dt} \int v^{3} w - \int F v_{t} - c_{1} \int F \frac{x}{|x|} \cdot \nabla v \\ - \frac{c_{2}}{8R^{3}} \int_{|y| \leq 2R} \int F \frac{(x - y)}{|x - y|} \cdot \nabla v  - \frac{c_{3}}{R} \int F \chi(\frac{|x|}{R}) x \cdot \nabla v \\ \lesssim \delta(\int \frac{1}{|x|} v^{4}) + \frac{1}{\delta} E(t) \| w \|_{L^{6}}^{2}.
\endaligned
\end{equation}
Next, splitting $F = -3v^{2} w - 3v w^{2}$, the Sobolev embedding theorem implies that
\begin{equation}\label{3.11}
-3  \int v_{t} v w^{2} dx \lesssim \| w \|_{L_{x}^{6}(\mathbf{R}^{3})}^{2} \| v \|_{L_{x}^{6}(\mathbf{R}^{3})} \| v_{t} \|_{L_{x}^{2}(\mathbf{R}^{3})} \lesssim E(t) \| w(t) \|_{L_{x}^{6}(\mathbf{R}^{3})}^{2}.
\end{equation}
Therefore,
\begin{equation}\label{3.26.5}
\aligned
\frac{d}{dt} \mathcal E(t) + c_{1} \pi v(t,0)^{2} + \frac{c_{2} \pi}{16 R^{3}} \int_{|y| \leq 2R} v(t,y)^{2} \\ + \frac{c_{1}}{2} \int \frac{1}{|x|} v^{4}  + \frac{c_{2}}{8 R^{3}}  \int_{|y| \leq 2R} \int \frac{1}{|x - y|} v^{4} \\ +\frac{c_{3}}{2R} \int \phi(\frac{|x|}{R}) [v_{t}^{2} + |\nabla v|^{2}]  + \frac{c_{3}}{4R} \int \chi(\frac{|x|}{R}) v^{4} \\
- \frac{d}{dt} \int v^{3} w + 3 \int v^{2} w v_{t} - c_{1} \int F \frac{x}{|x|} \cdot \nabla v \\ - \frac{c_{2}}{8R^{3}} \int_{|y| \leq 2R} \int F \frac{(x - y)}{|x - y|} \cdot \nabla v  - \frac{c_{3}}{R} \int F \chi(\frac{|x|}{R}) x \cdot \nabla v \\ \lesssim \delta(\int \frac{1}{|x|} v^{4}) + \frac{1}{\delta} E(t) \| w \|_{L^{6}}^{2}.
\endaligned
\end{equation}
Analysis of the other terms involving $-3 v^{2} w$ is similar.
\begin{equation}\label{3.25}
\int v w^{2} \frac{x}{|x|} \cdot \nabla v \lesssim \| w \|_{L^{6}}^{2} \| \nabla v \|_{L^{2}} \| v \|_{L^{6}} \lesssim E(t) \| w \|_{L^{6}}^{2},
\end{equation}
and
\begin{equation}\label{3.25.1}
 \frac{1}{8R^{3}} \int_{|y| \leq 2R} \int F \frac{(x - y)}{|x - y|} \cdot \nabla v \lesssim E(t) \| w \|_{L^{6}}^{2}.
 \end{equation}
 Since $\chi(\frac{|x|}{R}) \frac{x}{R}$ is also uniformly bounded,
 \begin{equation}\label{3.25.2}
 \frac{1}{R} \int v w^{2} \chi(\frac{|x|}{R}) x \cdot \nabla v \lesssim E(t) \| w \|_{L^{6}}^{2}.
 \end{equation}
 Therefore,
 \begin{equation}\label{3.25.3}
\aligned
\frac{d}{dt} \mathcal E(t) + c_{1} \pi v(t,0)^{2} + \frac{c_{2} \pi}{16 R^{3}} \int_{|y| \leq 2R} v(t,y)^{2} \\ + \frac{c_{1}}{2} \int \frac{1}{|x|} v^{4}  + \frac{c_{2}}{8 R^{3}}  \int_{|y| \leq 2R} \int \frac{1}{|x - y|} v^{4} \\ +\frac{c_{3}}{2R} \int \phi(\frac{|x|}{R}) [v_{t}^{2} + |\nabla v|^{2}]  + \frac{c_{3}}{4R} \int \chi(\frac{|x|}{R}) v^{4} \\
- \frac{d}{dt} \int v^{3} w + 3 \int v^{2} w v_{t} + 3c_{1} \int v^{2} w \frac{x}{|x|} \cdot \nabla v \\ + \frac{3 c_{2}}{8R^{3}} \int_{|y| \leq 2R} \int v^{2} w \frac{(x - y)}{|x - y|} \cdot \nabla v  + \frac{3 c_{3}}{R} \int v^{2} w \chi(\frac{|x|}{R}) x \cdot \nabla v \\ \lesssim \delta(\int \frac{1}{|x|} v^{4}) + \frac{1}{\delta} E(t) \| w \|_{L^{6}}^{2}.
\endaligned
\end{equation}

Next, by the product rule,
\begin{equation}\label{3.12}
3 \int v_{t} v^{2} w dx - \frac{d}{dt} \int v^{3} w dx = -\int v^{3} \partial_{t} w dx.
\end{equation}
Making a Littlewood--Paley decomposition,
\begin{equation}\label{3.13}
\aligned
\int v^{3} w_{t} dx = \sum_{j} \int v^{3} \partial_{t} w_{j} dx. 
\endaligned
\end{equation}
By Fourier support properties,
\begin{equation}\label{3.13.1}
\aligned
\int v^{3} \partial_{t} w_{j} dx = \int (v^{3} - (P_{\leq j - 3} v)^{3}) (\partial_{t} w_{j}) dx \\ = \int (P_{\geq j - 3} v)^{3} (\partial_{t} w_{j}) dx + 3 \int (P_{\geq j - 3} v)(P_{\leq j - 3} v) v \cdot \partial_{t} w_{j} dx.
\endaligned
\end{equation}
Using Lemma $\ref{l3.1}$,
\begin{equation}\label{3.17}
\aligned
\sum_{j} \int_{|x| \geq \frac{R}{2}} [v^{3} - (P_{\leq j - 3} v)^{3}] (\partial_{t} w_{j}) dx \\ \lesssim \sum_{j} (\| \frac{1}{|x|^{1/4}} |P_{\leq j} v| \|_{L^{4}} + \| \frac{1}{|x|^{1/4}} |P_{\geq j} v| \|_{L^{4}}) \| P_{\geq j - 3} v \|_{L_{x}^{2}} \| |x|^{1/2} \partial_{t} w_{j} \|_{L_{x}^{\infty}(|x| \geq \frac{R}{2})} \\ \lesssim (\int \frac{1}{|x|} v^{4})^{1/2} \sum_{j} \| P_{\geq j - 3} v \|_{L_{x}^{2}} \| |x|^{1/2} \partial_{t} w_{j} \|_{L_{x}^{\infty}(|x| \geq \frac{R}{2})}.
\endaligned
\end{equation}
By the Cauchy--Schwartz inequality,
\begin{equation}\label{3.15}
\aligned
(\ref{3.17}) \lesssim \delta(\int \frac{1}{|x|} v^{4}) + \frac{1}{\delta} (\sum_{j} \| P_{\geq j - 3} v \|_{L^{2}} \| \partial_{t} w_{j} \|_{L^{\infty}(|x| \geq \frac{R}{2})})^{2}.
\endaligned
\end{equation}
By Bernstein's inequality and Young's inequality,
\begin{equation}\label{3.21}
\aligned
(\sum_{j} \| P_{\geq j - 3} v \|_{L^{2}} \| |x|^{1/2} \partial_{t} w_{j} \|_{L^{\infty}(|x| \geq \frac{R}{2})})^{2} \\ \leq (\sum_{j} (\sum_{k \geq j - 3} 2^{k} 2^{j - k} \| P_{k} v \|_{L^{2}}) \cdot 2^{-j} \| |x|^{1/2} \partial_{t} w_{j} \|_{L^{\infty}(|x| \geq \frac{R}{2})})^{2} \\
\lesssim (\sum_{k} 2^{2k} \| P_{k} v \|_{L^{2}}^{2}) (\sum_{j} 2^{-2j} \| |x|^{1/2} \partial_{t} w_{j} \|_{L^{\infty}(|x| \geq \frac{R}{2})}^{2}) \\
\lesssim E(t) (\sum_{j} 2^{-2j} \| |x|^{1/2} \partial_{t} w_{j} \|_{L^{\infty}(|x| \geq \frac{R}{2})}^{2}).
\endaligned
\end{equation}
Therefore,
\begin{equation}\label{3.22}
\aligned
\frac{d}{dt} \mathcal E(t)  + c_{1} \pi v(t,0)^{2} + \frac{c_{2} \pi}{16 R^{3}} \int_{|y| \leq 2R} v(t,y)^{2} \\ + \frac{c_{1}}{2} \int \frac{1}{|x|} v^{4}  + \frac{c_{2}}{8 R^{3}}  \int_{|y| \leq 2R} \int \frac{1}{|x - y|} v^{4} \\ +\frac{c_{3}}{2R} \int \phi(\frac{|x|}{R}) [v_{t}^{2} + |\nabla v|^{2}]  + \frac{c_{3}}{4R} \int \chi(\frac{|x|}{R}) v^{4} \\
- \sum_{j} \int_{|x| \leq \frac{R}{2}} (v^{3} - (P_{\leq j - 3} v)^{3}) \cdot \partial_{t} w_{j} + 3c_{1} \int v^{2} w \frac{x}{|x|} \cdot \nabla v \\ + \frac{3 c_{2}}{8R^{3}} \int_{|y| \leq 2R} \int v^{2} w \frac{(x - y)}{|x - y|} \cdot \nabla v  + \frac{3 c_{3}}{R} \int v^{2} w \chi(\frac{|x|}{R}) x \cdot \nabla v \\ \lesssim \delta(\int \frac{1}{|x|} v^{4}) + \frac{1}{\delta} E(t) \| w \|_{L^{6}}^{2} + \frac{1}{\delta} E(t) (\sum_{j} 2^{-2j} \| |x|^{1/2} \partial_{t} w_{j} \|_{L^{\infty}(|x| \geq \frac{R}{2})}^{2}).
\endaligned
\end{equation}
By $(\ref{2.31})$, H{\"o}lder's inequality, and the Cauchy--Schwartz inequality,
\begin{equation}\label{3.24}
\aligned
\sum_{j} \int_{|x| \leq \frac{R}{2}} v (P_{\leq j - 3} v)(P_{\geq j - 3} v) \cdot \partial_{t} w_{j} dx \\ \leq \sum_{j} \| \partial_{t} w_{j} \|_{L^{6}} \| v \|_{L^{6}(|x| \leq \frac{R}{2})} \| P_{\geq j - 3} v \|_{L^{2}} \| P_{\leq j - 3} v \|_{L^{6}} \\ \lesssim \delta R E(t) (\frac{1}{R} \int_{|x| \leq R} |\nabla v|^{2} + \frac{1}{R^{3}} \int_{|x| \leq R} v^{2}) + \frac{1}{\delta} (\sum_{j} \| P_{\geq j - 3} v \|_{L^{2}} \| \partial_{t} w_{j} \|_{L^{6}})^{2}.
\endaligned
\end{equation}
Following $(\ref{3.21})$,
\begin{equation}\label{3.24.1}
(\sum_{j} \| P_{\geq j - 3} v \|_{L^{2}} \| \partial_{t} w_{j} \|_{L^{6}})^{2} \lesssim E(t) (\sum_{j} 2^{-2j} \| \partial_{t} w_{j} \|_{L^{6}}^{2}).
\end{equation}

Next, following $(\ref{3.21})$, by the Cauchy--Schwartz inequality, and Lemma $\ref{l2.2}$,
\begin{equation}\label{3.28}
\aligned
\sum_{j} \int_{|x| \leq \frac{R}{2}} (P_{\geq j - 3} v)^{3} \cdot \partial_{t} w_{j} dx \\ \lesssim \sum_{j} \| P_{\geq j - 3} v \|_{L^{4}(|x| \leq \frac{R}{2})}^{2} \| P_{\geq j - 3} v \|_{L^{3}(\mathbf{R}^{3})} \| \partial_{t} w_{j} \|_{L^{6}(\mathbf{R}^{3})} \\
\lesssim \delta R \| \nabla v(t) \|_{L^{2}}^{2} [\frac{1}{R} \| \nabla v \|_{L^{2}(|x| \leq R)}^{2} + \frac{1}{R^{3}} \| v \|_{L^{2}(|x| \leq R)}^{2}] + \delta (\int \frac{1}{|x|} v^{4}) \\ + \frac{1}{\delta} (\sum_{j} 2^{-j/2} \| P_{\geq j - 3} v \|_{L^{3}} \| \partial_{t} w_{j} \|_{L^{6}})^{2} \\
\lesssim \delta R E(t) [\frac{1}{R} \| \nabla v \|_{L^{2}(|x| \leq R)}^{2} + \frac{1}{R^{3}} \| v \|_{L^{2}(|x| \leq R)}^{2}] \\ + \delta (\int \frac{1}{|x|} v^{4}) + \frac{1}{\delta} E(t) [\sum_{j} 2^{-2j} \| \partial_{t} w_{j} \|_{L^{6}}^{2}].
\endaligned
\end{equation}
Therefore,
\begin{equation}\label{3.26}
\aligned
\frac{d}{dt} \mathcal E(t)  + c_{1} \pi v(t,0)^{2} + \frac{c_{2} \pi}{16 R^{3}} \int_{|y| \leq 2R} v(t,y)^{2} \\ + \frac{c_{1}}{2} \int \frac{1}{|x|} v^{4}  + \frac{c_{2}}{8 R^{3}}  \int_{|y| \leq 2R} \int \frac{1}{|x - y|} v^{4} \\ +\frac{c_{3}}{2R} \int \phi(\frac{|x|}{R}) [v_{t}^{2} + |\nabla v|^{2}]  + \frac{c_{3}}{4R} \int \chi(\frac{|x|}{R}) v^{4} \\
+ 3c_{1} \int v^{2} w \frac{x}{|x|} \cdot \nabla v \\ + \frac{3 c_{2}}{8R^{3}} \int_{|y| \leq 2R} \int v^{2} w \frac{(x - y)}{|x - y|} \cdot \nabla v  + \frac{3 c_{3}}{R} \int v^{2} w \chi(\frac{|x|}{R}) x \cdot \nabla v \\ \lesssim \delta(\int \frac{1}{|x|} v^{4}) + \delta R E(t) [\frac{1}{R} \| \nabla v \|_{L^{2}(|x| \leq R)}^{2} + \frac{1}{R^{3}} \| v \|_{L^{2}(|x| \leq R)}^{2}] \\ + \frac{1}{\delta} E(t) \| w \|_{L^{6}}^{2} + \frac{1}{\delta} E(t)(\sum_{j} 2^{-2j} \| \partial_{t} w_{j} \|_{L^{6}}^{2}) \\ + \frac{1}{\delta} E(t) (\sum_{j} 2^{-2j} \| |x|^{1/2} \partial_{t} w_{j} \|_{L^{\infty}(|x| \geq \frac{R}{2})}^{2}).
\endaligned
\end{equation}
Integrating by parts,
\begin{equation}\label{3.30}
3c_{1} \int v^{2} w \frac{x}{|x|} \cdot \nabla v dx = -2c_{1} \int \frac{1}{|x|} v^{3} w - c_{1} \int v^{3} (\nabla w) \cdot \frac{x}{|x|}.
\end{equation}
Following $(\ref{3.26.2})$,
\begin{equation}\label{3.31}
-2c_{1} \int \frac{1}{|x|} v^{3} w dx \lesssim \frac{1}{\delta} E(t) \| w(t) \|_{L_{x}^{6}}^{2} + \delta (\int \frac{1}{|x|} v^{4} dx).
\end{equation}
The term
\begin{equation}
-c_{1} \int (v^{3} - (P_{\leq j - 3} v)^{3}) (\nabla w_{j}) \cdot \frac{x}{|x|} dx
\end{equation}
may be estimated using exactly the same arguments as in the estimates for $(\ref{3.17})$.

Now, the Fourier support of $(\nabla w_{j})(P_{\leq j - 3} v)^{3}$ is $|\xi| \sim 2^{j}$, so integrating by parts,
\begin{equation}\label{3.32}
\aligned
c \int (P_{\leq j - 3} v)^{3} (\nabla w_{j}) \cdot \frac{x}{|x|} dx = \int \frac{x_{l} x_{k}}{|x|^{3}} \frac{\partial_{k}}{\Delta} (P_{\leq j - 3} v)^{3} (\partial_{l} w_{j}) \\
\lesssim 2^{-j} \| \frac{1}{|x|^{1/4}} P_{\leq j - 3} v \|_{L^{4}}^{2} \| \frac{1}{|x|^{1/2}} P_{\leq j - 3} v \|_{L^{10/3}} \| \partial_{k} w_{j} \|_{L^{5}}.
\endaligned
\end{equation}
Then by the Cauchy--Schwartz inequality,
\begin{equation}\label{3.33}
\aligned
\sum_{j} c \int (P_{\leq j - 3} v)^{3} (\nabla w_{j}) \cdot \frac{x}{|x|} dx \\ \lesssim \delta (\int \frac{1}{|x|} v^{4}) + \frac{1}{\delta} (\sum_{j} 2^{-j} \| \frac{1}{|x|^{1/2}} P_{\leq j - 3} v \|_{L^{10/3}} \| \nabla w_{j} \|_{L^{5}})^{2},
\endaligned
\end{equation}
and then by Bernstein's inequality,
\begin{equation}\label{3.34}
\lesssim \delta (\int \frac{1}{|x|} v^{4}) + \frac{1}{\delta} E(t) (\sum_{j} 2^{j/5} \| w_{j} \|_{L^{5}}^{2}).
\end{equation}
Therefore,
\begin{equation}\label{3.35}
\aligned
\frac{d}{dt} \mathcal E(t)  + c_{1} \pi v(t,0)^{2} + \frac{c_{2} \pi}{16 R^{3}} \int_{|y| \leq 2R} v(t,y)^{2} \\ + \frac{c_{1}}{2} \int \frac{1}{|x|} v^{4}  + \frac{c_{2}}{8 R^{3}}  \int_{|y| \leq 2R} \int \frac{1}{|x - y|} v^{4} \\ +\frac{c_{3}}{2R} \int \phi(\frac{|x|}{R}) [v_{t}^{2} + |\nabla v|^{2}]  + \frac{c_{3}}{4R} \int \chi(\frac{|x|}{R}) v^{4} \\
 + \frac{3 c_{2}}{8R^{3}} \int_{|y| \leq 2R} \int v^{2} w \frac{(x - y)}{|x - y|} \cdot \nabla v  + \frac{3 c_{3}}{R} \int v^{2} w \chi(\frac{|x|}{R}) x \cdot \nabla v \\ \lesssim \delta(\int \frac{1}{|x|} v^{4}) + \delta R E(t) [\frac{1}{R} \| \nabla v \|_{L^{2}(|x| \leq R)}^{2} + \frac{1}{R^{3}} \| v \|_{L^{2}(|x| \leq R)}^{2}] \\ + \frac{1}{\delta} E(t) \| w \|_{L^{6}}^{2} + \frac{1}{\delta} E(t) (\sum_{j} 2^{j/5} \| w_{j} \|_{L^{5}}^{2})+ \frac{1}{\delta} E(t)(\sum_{j} 2^{-2j} \| \partial_{t} w_{j} \|_{L^{6}}^{2}) \\ + \frac{1}{\delta} E(t) (\sum_{j} 2^{-2j} \| |x|^{1/2} \partial_{t} w_{j} \|_{L^{\infty}(|x| \geq \frac{R}{2})}^{2}) \\ 
\frac{1}{\delta} E(t)(\sum_{j} 2^{-2j} \| \nabla w_{j} \|_{L^{6}}^{2}) + \frac{1}{\delta} E(t) (\sum_{j} 2^{-2j} \| |x|^{1/2} \nabla w_{j} \|_{L^{\infty}(|x| \geq \frac{R}{2})}^{2}) .
\endaligned
\end{equation}
Like $\frac{x}{|x|}$ the potentials
\begin{equation}\label{3.36}
a(x) = \chi(\frac{2 |x|}{R}) \frac{x}{R}, \qquad \text{and} \qquad a(x) = \int_{|y| \leq 2R} \frac{(x - y)}{|x - y|} dy
\end{equation}
are also bounded, radial functions satisfying
\begin{equation}\label{3.37}
\nabla \cdot a(x) \lesssim \frac{1}{|x|},
\end{equation}
and therefore, the analysis of 
\begin{equation}\label{3.38}
+ \frac{3 c_{2}}{8R^{3}} \int_{|y| \leq 2R} \int v^{2} w \frac{(x - y)}{|x - y|} \cdot \nabla v  + \frac{3 c_{3}}{R} \int v^{2} w \chi(\frac{|x|}{R}) x \cdot \nabla v
\end{equation}
may be carried out in much the same manner as
\begin{equation}\label{3.39}
\int v^{2} w \frac{x}{|x|} \cdot \nabla v.
\end{equation}
Therefore,
\begin{equation}\label{3.40}
\aligned
\frac{d}{dt} \mathcal E(t) 
 \lesssim  \frac{1}{\delta} E(t) \| w \|_{L^{6}}^{2} + \frac{1}{\delta} E(t) (\sum_{j} 2^{j/5} \| w_{j} \|_{L^{5}}^{2})+ \frac{1}{\delta} E(t)(\sum_{j} 2^{-2j} \| \partial_{t} w_{j} \|_{L^{6}}^{2}) \\ + \frac{1}{\delta} E(t) (\sum_{j} 2^{-2j} \| |x|^{1/2} \partial_{t} w_{j} \|_{L^{\infty}(|x| \geq \frac{R}{2})}^{2}) \\ 
\frac{1}{\delta} E(t)(\sum_{j} 2^{-2j} \| \nabla w_{j} \|_{L^{6}}^{2}) + \frac{1}{\delta} E(t) (\sum_{j} 2^{-2j} \| |x|^{1/2} \nabla w_{j} \|_{L^{\infty}(|x| \geq \frac{R}{2})}^{2}) .
\endaligned
\end{equation}
Since $E(t) \sim \mathcal E(t)$,
\begin{equation}\label{3.41}
\aligned
\frac{d}{dt} \ln(\mathcal E(t))
 \lesssim  \frac{1}{\delta} \| w \|_{L^{6}}^{2} + \frac{1}{\delta} (\sum_{j} 2^{j/5} \| w_{j} \|_{L^{5}}^{2})+ \frac{1}{\delta} (\sum_{j} 2^{-2j} \| \partial_{t} w_{j} \|_{L^{6}}^{2}) \\ + \frac{1}{\delta} (\sum_{j} 2^{-2j} \| |x|^{1/2} \partial_{t} w_{j} \|_{L^{\infty}(|x| \geq \frac{1}{2 \mathcal E(T)})}^{2}) \\ 
\frac{1}{\delta} (\sum_{j} 2^{-2j} \| \nabla w_{j} \|_{L^{6}}^{2}) + \frac{1}{\delta} (\sum_{j} 2^{-2j} \| |x|^{1/2} \nabla w_{j} \|_{L^{\infty}(|x| \geq \frac{1}{2 \mathcal E(T)})}^{2}) .
\endaligned
\end{equation}
Integrating in time and combining $(\ref{2.18.1})$ with $(\ref{3.4.1})$,
\begin{equation}\label{3.45}
\aligned
\ln(\mathcal E(T)) - \ln(\mathcal E(0)) \lesssim \frac{\epsilon^{2}}{\delta} \ln(T) + \frac{\epsilon^{2}}{\delta} \ln(\mathcal E(T)) + \epsilon.
\endaligned
\end{equation}
Doing some algebra,
\begin{equation}\label{3.46}
\aligned
\ln(\mathcal E(T)) \leq (\frac{1}{1 - \frac{C \epsilon^{2}}{\delta}}) \ln(\mathcal E(0)) + \frac{C \epsilon^{2}}{\delta(1 - \frac{C \epsilon^{2}}{\delta})} \ln(T) + \frac{C \epsilon}{(1 - \frac{C \epsilon^{2}}{\delta})}.
\endaligned
\end{equation}
This proves that for any $t$, $E(t) \sim \mathcal E(t) \lesssim (1 + t)^{C \epsilon}$. $\Box$

\section{Proof of scattering}
By time reversal symmetry, it suffices to prove
\begin{theorem}\label{t4.0}
For any radial initial data $(u_{0}, u_{1}) \in \dot{H}^{1/2} \times \dot{H}^{-1/2}$, the solution to $(\ref{1.1})$ scatters forward in time.
\end{theorem}
This theorem is proved using hyperbolic coordinates. By the dominated convergence theorem, there exists $R(\epsilon) < \infty$ such that
\begin{equation}\label{4.1}
\| S(t)(u_{0}, u_{1}) \|_{L_{t,x}^{4}(|x| \geq R + |t|)} < \epsilon.
\end{equation}
Then by finite propagation speed and small data arguments, if $u$ is a global solution to $(\ref{1.1})$, then
\begin{equation}\label{4.2}
\| u \|_{L_{t,x}^{4}(|x| \geq R + |t|)} \lesssim \epsilon.
\end{equation}
Rescaling, $(u_{0}(x), u_{1}(x)) \mapsto (2 R u_{0}(2 R x), (2 R)^{2} u_{1}(2R x))$,
\begin{equation}\label{4.3}
\| u \|_{L_{t,x}^{4}(|x| \geq \frac{1}{2} + |t|)} \lesssim \epsilon.
\end{equation}

The quantity
\begin{equation}\label{4.4}
\| u \|_{L_{t,x}^{4}([0, \infty) \times \{|x| \leq \frac{1}{2} + t \})}
\end{equation}
is estimated using hyperbolic coordinates, which combined with $(\ref{4.3})$ proves
\begin{equation}\label{4.5}
\| u \|_{L_{t,x}^{4}([0, \infty) \times \mathbf{R}^{3})} < \infty.
\end{equation}

Make a time translation so that
\begin{equation}\label{4.6}
u(1, x) = 2 R u_{0}(2 R x), \qquad u_{t}(1,x) = (2 R)^{2} u_{1}(2 Rx).
\end{equation}
After time translation, $(\ref{4.3})$ implies
\begin{equation}\label{4.7}
\| u \|_{L_{t,x}^{4}([1, \infty) \times \{ |x| \geq t - \frac{1}{2} \})} \lesssim \epsilon.
\end{equation}
Switching to hyperbolic coordinates for the region inside the cone, let
\begin{equation}\label{4.8}
\tilde{u}(\tau, s) = \frac{e^{\tau} \sinh s}{s} u(e^{\tau} \cosh s, e^{\tau} \sinh s).
\end{equation}
Then making a change of variables,
\begin{equation}\label{4.9}
\aligned
\int_{0}^{\infty} \int_{0}^{\infty} \tilde{u}(\tau, s)^{4} (\frac{s}{\sinh s})^{2} s^{2} ds d\tau \\ = \int_{0}^{\infty} \int_{0}^{\infty} u(e^{\tau} \cosh s, e^{\tau} \sinh s)^{4} e^{2 \tau} \sinh^{2} e^{2 \tau} ds d\tau \\ = \int_{1}^{\infty} \int_{t^{2} - r^{2} \geq 1} u(t,r)^{4} r^{2} dr dt \geq \int_{1}^{\infty} \int_{t \geq r} u(t,r)^{4} r^{2} dr dt.
\endaligned
\end{equation}
Therefore,
\begin{equation}\label{4.10}
\int_{0}^{\infty} \int_{0}^{\infty} \tilde{u}(\tau, s)^{4} (\frac{s}{\sinh s})^{2} s^{2} ds d\tau < \infty,
\end{equation}
combined with $(\ref{4.7})$ implies
\begin{equation}\label{4.11}
\| u \|_{L_{t,x}^{4}([1, \infty) \times \mathbf{R}^{3})} < \infty,
\end{equation}
which after undoing time translation, implies $(\ref{4.5})$. Also, by direct computation,
\begin{equation}\label{4.12}
(\partial_{\tau \tau} - \partial_{ss} - \frac{2}{s} \partial_{s}) \tilde{u}(\tau, s) + (\frac{s}{\sinh s})^{2} \tilde{u}^{3} = 0,
\end{equation}
with
\begin{equation}\label{4.13}
\tilde{u}|_{\tau = 0} = \frac{e^{\tau} \sinh s}{s} u(e^{\tau} \cosh s, e^{\tau} \sinh s)|_{\tau = 0},
\end{equation}
and
\begin{equation}\label{4.14}
\tilde{u}_{\tau}|_{\tau = 0} = \partial_{\tau} (\frac{e^{\tau} \sinh s}{s} u(e^{\tau} \cosh s, e^{\tau} \sinh s))|_{\tau = 0}.
\end{equation}
A solution to $(\ref{4.12})$ has the conserved energy,
\begin{equation}\label{4.15}
E(\tau) = \frac{1}{2} \| \tilde{u}_{\tau} \|_{L^{2}}^{2} + \frac{1}{2} \| \tilde{u}_{s} \|_{L^{2}}^{2} + \frac{1}{4} \int \tilde{u}(\tau, s)^{4} (\frac{s}{\sinh s})^{2} s^{2} ds.
\end{equation}
For now, assume the following lemma.
\begin{lemma}\label{l4.1}
There exists a decomposition
\begin{equation}\label{4.16}
\tilde{u}|_{\tau = 0} = \frac{e^{\tau} \sinh s}{s} u(e^{\tau} \cosh s, e^{\tau} \sinh s)|_{\tau = 0} = \tilde{v}_{0} + \tilde{w}_{0}
\end{equation}
and
\begin{equation}\label{4.17}
\tilde{u}_{\tau}|_{\tau = 0} = \partial_{\tau} (\frac{e^{\tau} \sinh s}{s} u(e^{\tau} \cosh s, e^{\tau} \sinh s))|_{\tau = 0} = \tilde{v}_{1} + \tilde{w}_{1},
\end{equation}
with
\begin{equation}\label{4.18}
\frac{1}{2} \int |\partial_{s} \tilde{v}_{0}|^{2} s^{2} + \frac{1}{2} \int |\tilde{v}_{1}|^{2} s^{2} + \frac{1}{4} \int \tilde{v}_{0}^{4} (\frac{s}{\sinh s})^{2} s^{2} < \infty,
\end{equation}
and
\begin{equation}\label{4.19}
\| \tilde{w}_{0} \|_{\dot{H}^{1/2}} + \| w_{1} \|_{\dot{H}^{-1/2}} \leq \epsilon.
\end{equation}
\end{lemma}
\noindent \textbf{Remark:} Following $(\ref{3.1})$, it is enough to prove $\tilde{u}_{0} \in \dot{H}^{1} + \dot{H}^{1/2}$ and $\tilde{u}_{1} \in L^{2} + \dot{H}^{-1/2}$ and then truncate in frequency.\vspace{5mm}

\noindent \emph{Proof of Theorem $\ref{t4.0}$:} Let $\tilde{v}$ and $\tilde{w}$ solve
\begin{equation}\label{4.19.2}
(\partial_{\tau \tau} - \Delta) \tilde{w} + (\frac{s}{\sinh s})^{2} \tilde{w}^{3} = 0, \qquad \tilde{w}(0,y) = \tilde{w}_{0}, \qquad \tilde{w}_{\tau}(0,y) = \tilde{w}_{1},
\end{equation}
and
\begin{equation}\label{4.19.1}
(\partial_{\tau \tau} - \Delta) \tilde{v} + (\frac{s}{\sinh s})^{2} [\tilde{v}^{3} + 3 \tilde{v}^{2} \tilde{w} + 3 \tilde{v} \tilde{w}^{2}] = 0, \qquad \tilde{v}(0,y) = \tilde{v}_{0}, \qquad \tilde{v}_{\tau}(0,y) = \tilde{v}_{1}.
\end{equation}
Define the energy,
\begin{equation}\label{4.18}
E(\tau) = \frac{1}{2} \int |\partial_{s} \tilde{v}|^{2} s^{2} + \frac{1}{2} \int |\partial_{\tau} \tilde{v}|^{2} s^{2} + \frac{1}{4} \int \tilde{v}^{4} (\frac{s}{\sinh s})^{2} s^{2}.
\end{equation}
As in the proof of global well-posedness, define the quantity
\begin{equation}\label{4.19}
\mathcal E(\tau) = E(\tau) + M(\tau) + \int \tilde{v}^{3} \tilde{w} (\frac{s}{\sinh s})^{2} s^{2} ds,
\end{equation}
where
\begin{equation}\label{4.20}
M(\tau) = c \int \tilde{v}_{\tau} \tilde{v}_{s} s^{2} ds + c \int \tilde{v}_{\tau}  \tilde{v} ds.
\end{equation}
Then by direct computation, making a slight modification of $(\ref{2.24})$ and $(\ref{2.36})$,
\begin{equation}\label{4.21}
\aligned
\frac{d}{d\tau} M(\tau) = -\frac{1}{2} \tilde{v}(\tau, 0)^{2} - \frac{1}{2} \int \tilde{v}(\tau, s)^{4} (\frac{s}{\sinh s})^{2} (\frac{\cosh s}{\sinh s}) s^{2} ds \\
-3 \int \tilde{v}^{2} \tilde{v}_{s} \tilde{w} (\frac{s}{\sinh s})^{2} s^{2} ds - 3 \int \tilde{v} \tilde{v}_{s} \tilde{w}^{2} (\frac{s}{\sinh s})^{2} s^{2} ds \\
- 3 \int \tilde{v}^{3} \tilde{w} (\frac{s}{\sinh s})^{2} s ds - 3 \int \tilde{v}^{2} \tilde{w}^{2} (\frac{s}{\sinh s})^{2} s ds.
\endaligned
\end{equation}
Therefore,
\begin{equation}\label{4.22}
\aligned
\frac{d}{d\tau} \mathcal E(\tau) = -\frac{c}{2} \tilde{v}(\tau, 0)^{2} - \frac{c}{2} \int \tilde{v}(\tau, s)^{4} (\frac{s}{\sinh s})^{2} (\frac{\cosh s}{\sinh s}) s^{2} ds \\
-3c \int \tilde{v}^{2} \tilde{v}_{s} \tilde{w} (\frac{s}{\sinh s})^{2} s^{2} ds - 3c \int \tilde{v} \tilde{v}_{s} \tilde{w}^{2} (\frac{s}{\sinh s})^{2} s^{2} ds \\
- 3c \int \tilde{v}^{3} \tilde{w} (\frac{s}{\sinh s})^{2} s ds - 3c \int \tilde{v}^{2} \tilde{w}^{2} (\frac{s}{\sinh s})^{2} s ds \\
-3 \int (\frac{s}{\sinh s})^{2} \tilde{v}^{3} \tilde{w}_{\tau} s ds - 3 \int (\frac{s}{\sinh s})^{2} \tilde{v} \tilde{v}_{\tau} \tilde{w}^{2} s ds.
\endaligned
\end{equation}
By Hardy's inequality,
\begin{equation}\label{4.23}
\aligned
-3c \int \tilde{v} \tilde{v}_{s} \tilde{w}^{2} (\frac{s}{\sinh s})^{2} s^{2} ds - 3c \int \tilde{v}^{2} \tilde{w}^{2} (\frac{s}{\sinh s})^{2} s ds\\ - 3 \int (\frac{s}{\sinh s})^{2} \tilde{v} \tilde{v}_{\tau} \tilde{w}^{2} s^{2} ds \lesssim E(\tau) \| \tilde{w} \|_{L^{6}}^{2}.
\endaligned
\end{equation}
Also, by Hardy's inequality and the Cauchy--Schwartz inequality,
\begin{equation}\label{4.24.2}
\aligned
\int \tilde{v}^{3} \tilde{w} (\frac{s}{\sinh s})^{2} s ds \lesssim \delta(\int \tilde{v}^{4} (\frac{\cosh s}{\sinh s}) (\frac{s}{\sinh s})^{2} s^{2} ds) + \frac{1}{\delta} \| \tilde{w} \|_{L^{6}}^{2} \| \frac{1}{|x|^{1/2}} \tilde{v} \|_{L^{3}}^{2} \\ \lesssim \delta(\int (\frac{s}{\sinh s})^{2} \frac{\cosh s}{\sinh s} \tilde{v}^{4} s^{2} ds) + \frac{1}{\delta} \| \tilde{w} \|_{L^{6}}^{2} E(\tau).
\endaligned
\end{equation}
Therefore,
\begin{equation}\label{4.24}
\aligned
\frac{d}{d\tau} \mathcal E(\tau) + \frac{c}{2} \tilde{v}(\tau, 0)^{2} + \frac{c}{2} \int \tilde{v}(\tau, s)^{4} (\frac{s}{\sinh s})^{2} (\frac{\cosh s}{\sinh s}) s^{2} ds \\
+ 3c \int \tilde{v}^{2} \tilde{w} (\frac{s}{\sinh s})^{2} \tilde{v}_{s} s^{2} ds + 3 \int (\frac{s}{\sinh s})^{2} \tilde{v}^{3} \tilde{w}_{\tau} s^{2} ds \\ \lesssim \frac{1}{\delta} E(\tau) \| w \|_{L^{6}}^{2} + \delta(\int \tilde{v}^{4} (\frac{s}{\sinh s})^{2} (\frac{\cosh s}{\sinh s}) s^{2} ds).
\endaligned
\end{equation}
Integrating by parts,
\begin{equation}\label{4.25}
3c \int \tilde{v}^{2} \tilde{v}_{s} \tilde{w} (\frac{s}{\sinh s})^{2} s^{2} ds = -c \int \tilde{v}^{3} \tilde{w}_{s} (\frac{s}{\sinh s})^{2} s^{2} ds - c \int \tilde{v}^{3} \tilde{w} \cdot \partial_{s}(\frac{s^{4}}{(\sinh s)^{2}}) ds.
\end{equation}
\noindent Since
\begin{equation}\label{4.26}
\partial_{s}(\frac{s^{4}}{(\sinh s)^{2}}) \lesssim s,
\end{equation}
by $(\ref{4.24.2})$,
\begin{equation}\label{4.27}
c \int \tilde{v}^{3} \tilde{w} \cdot \partial_{s} (\frac{s^{4}}{(\sinh s)^{2}}) \lesssim \delta(\int \tilde{v}^{4} (\frac{s}{\sinh s})^{2} (\frac{\cosh s}{\sinh s}) s^{2} ds) + \frac{1}{\delta} \| \tilde{w} \|_{L^{6}}^{2} E(\tau).
\end{equation}
Next, following $(\ref{3.13.1})$--$(\ref{3.15})$ and using Lemma $\ref{l3.1}$,
\begin{equation}\label{4.28}
\aligned
-c \sum_{j} \int_{s \geq \frac{R}{2}} [\tilde{v}^{3} - (P_{\leq j - 3} \tilde{v})^{3}] (\partial_{s} \tilde{w}_{j}) \cdot (\frac{s}{\sinh s})^{2} s^{2} ds \\ + 3 \sum_{j} \int_{s \geq \frac{R}{2}}  [\tilde{v}^{3} - (P_{\leq j - 3} \tilde{v})^{3}] (\partial_{\tau} \tilde{w}_{j}) (\frac{s}{\sinh s})^{2} s^{2} ds \\
\lesssim \delta (\int (\frac{\cosh s}{\sinh s}) (\frac{s}{\sinh s})^{2} \tilde{v}^{4} s^{2} ds) \\ + \frac{1}{\delta} E(\tau) (\sum_{j} 2^{-2j} \| (\nabla_{\tau, x} \tilde{w}_{j}) (\frac{\sinh s}{\cosh s})^{1/2} (\frac{s}{\sinh s}) \|_{L^{\infty}(|x| \geq \frac{R}{2})}^{2}).
\endaligned
\end{equation}
Next, by H{\"o}lder's inequality,
\begin{equation}\label{4.29}
\aligned
\sum_{j} \| (\tilde{v}^{3} - (P_{\leq j - 3} \tilde{v})^{3}) (\nabla_{\tau, x} \tilde{w}_{j}) \|_{L^{1}(|x| \leq \frac{R}{2})} \\ \lesssim \sum_{j} \| \tilde{v} \|_{L^{\infty}} \| P_{\geq j - 3} \tilde{v} \|_{L^{2}} \| \nabla_{\tau, x} \tilde{w}_{j} \|_{L^{6}} \| \tilde{v} \|_{L^{3}(|x| \leq \frac{R}{2})} \\
\lesssim E(\tau) (\sum_{j} 2^{-2j} \| \nabla_{\tau, x} \tilde{w}_{j} \|_{L^{6}}^{2}) + R E(\tau) \| v \|_{L^{\infty}}^{2}.
\endaligned
\end{equation}
Following $(\ref{3.32})$ and $(\ref{3.33})$,
\begin{equation}\label{4.30}
\aligned
\int (P_{\leq j - 3} \tilde{v})^{3} (\partial_{s} \tilde{w}_{j}) \cdot (\frac{s}{\sinh s})^{2} s^{2} ds + \int (P_{\leq j - 3} \tilde{v})^{3} (\partial_{\tau} \tilde{w}_{j}) \cdot (\frac{s}{\sinh s})^{2} s^{2} ds \\
\lesssim \delta (\int \frac{1}{|x|} \tilde{v}^{4}) + \frac{1}{\delta} E(\tau) (\sum_{j} 2^{-8j/5} \| \nabla_{\tau, x} \tilde{w}_{j} \|_{L^{5}}^{2}).
\endaligned
\end{equation}
Therefore,
\begin{equation}\label{4.24.1}
\aligned
\frac{d}{d\tau} \mathcal E(\tau) + \frac{c}{2} \tilde{v}(\tau, 0)^{2} + \frac{c}{2} \int \tilde{v}(\tau, s)^{4} (\frac{s}{\sinh s})^{2} (\frac{\cosh s}{\sinh s}) s^{2} ds \\ \lesssim E(\tau) (\sum_{j} 2^{-2j} \| \nabla_{\tau, x} w_{j} \|_{L^{6}}^{2}) + \frac{1}{\delta} E(\tau) (\sum_{j} 2^{-8j/5} \| \nabla_{\tau, x} \tilde{w}_{j} \|_{L^{5}}^{2}) + R E(\tau) \| v \|_{L^{\infty}}^{2} \\
+ \frac{1}{\delta} E(\tau) \| w \|_{L^{6}}^{2} + \delta(\int \tilde{v}^{4} (\frac{s}{\sinh s})^{2} (\frac{\cosh s}{\sinh s}) s^{2} ds).
\endaligned
\end{equation}
Absorbing
\begin{equation}
\delta(\int \tilde{v}^{4} (\frac{s}{\sinh s})^{2} (\frac{\cosh s}{\sinh s}) s^{2} ds)
\end{equation}
into the left hand side,
\begin{equation}
\aligned
\frac{d}{d\tau} \mathcal E(\tau) + \frac{c}{4} \int \tilde{v}(\tau, s)^{4} (\frac{s}{\sinh s})^{2} (\frac{\cosh s}{\sinh s}) s^{2} ds \\ \lesssim E(\tau) (\sum_{j} 2^{-2j} \| \nabla_{\tau, x} w_{j} \|_{L^{6}}^{2}) + \frac{1}{\delta} E(\tau) (\sum_{j} 2^{-8j/5} \| \nabla_{\tau, x} \tilde{w}_{j} \|_{L^{5}}^{2}) \\ + R E(\tau) \| v \|_{L^{\infty}}^{2} + \frac{1}{\delta} E(\tau) \| w \|_{L^{6}}^{2}.
\endaligned
\end{equation}
Since $E(\tau) \sim \mathcal E(\tau)$,
\begin{equation}
\aligned
\frac{d}{d\tau} \ln(\mathcal E(\tau)) + \frac{c}{4 \mathcal E(\tau)} \int \tilde{v}(\tau, s)^{4} (\frac{s}{\sinh s})^{2} (\frac{\cosh s}{\sinh s}) s^{2} ds \\ \lesssim (\sum_{j} 2^{-2j} \| \nabla_{\tau, x} w_{j} \|_{L^{6}}^{2}) + \frac{1}{\delta} (\sum_{j} 2^{-8j/5} \| \nabla_{\tau, x} \tilde{w}_{j} \|_{L^{5}}^{2}) \\ + R \| v \|_{L^{\infty}}^{2} + \frac{1}{\delta} \| w \|_{L^{6}}^{2} + \frac{1}{\delta} (\sum_{j} 2^{-2j} \| (\nabla_{\tau, x} \tilde{w}_{j}) (\frac{\sinh s}{\cosh s})^{1/2} (\frac{s}{\sinh s}) \|_{L^{\infty}(s \geq \frac{R}{2})}^{2}).
\endaligned
\end{equation}
Suppose $T$ is such that $\mathcal E(T) = \sup_{0 < \tau < T} \mathcal E(\tau)$. Integrating in $\tau$,
\begin{equation}\label{4.28.1}
\aligned
\ln(\mathcal E(T)) - \ln(\mathcal E(0)) + \frac{c}{4} \int_{0}^{T} \frac{1}{\mathcal E(\tau)} \int \tilde{v}(\tau, s)^{4} (\frac{s}{\sinh s})^{2} (\frac{\cosh s}{\sinh s}) s^{2} ds d\tau \\ \lesssim \frac{\epsilon^{2}}{\delta} (1 - \ln(R)) + \epsilon^{2} + \int_{0}^{T} R \| \tilde{v} \|_{L^{\infty}}^{2} d\tau.
\endaligned
\end{equation}
Now by direct computation,
\begin{equation}\label{4.29}
\| (\frac{s}{\sinh s})^{1/2} \tilde{u} \|_{L^{4}} \lesssim \| (\frac{s}{\sinh s})^{1/2} (\frac{\cosh s}{\sinh s})^{1/4} \tilde{v} \|_{L^{4}} + \| \tilde{w} \|_{L^{4}}.
\end{equation}
If $I$ is an interval on which $\| (\frac{s}{\sinh s})^{1/2} \tilde{u} \|_{L_{\tau, x}^{4}(I)} \lesssim \epsilon$, then by $(\ref{6.10})$ and $(\ref{4.19.1})$,
\begin{equation}\label{4.30}
\| \tilde{v} \|_{L_{\tau}^{2} L_{x}^{\infty}(I \times \mathbf{R}^{3})} \lesssim \| \nabla \tilde{v} \|_{L_{\tau}^{\infty} L_{x}^{2}} + \| \tilde{v}_{\tau} \|_{L_{\tau}^{\infty} L_{x}^{2}} + \| \tilde{v} \|_{L_{\tau}^{2} L_{x}^{\infty}} (\int_{I} \int \tilde{u}^{4} (\frac{s}{\sinh s})^{2} s^{2} ds d\tau)^{1/2},
\end{equation}
which implies
\begin{equation}\label{4.30.1}
\| \tilde{v} \|_{L_{\tau}^{2} L_{x}^{\infty}(I \times \mathbf{R}^{3})} \lesssim \| \nabla \tilde{v} \|_{L_{\tau}^{\infty} L_{x}^{2}} + \| \tilde{v}_{\tau} \|_{L_{\tau}^{\infty} L_{x}^{2}},
\end{equation}
and therefore,
\begin{equation}\label{4.31}
\int_{0}^{T} R \| \tilde{v} \|_{L^{\infty}}^{2} d\tau \lesssim R \mathcal E(T) (\int_{0}^{T} \int \tilde{v}^{4} (\frac{s}{\sinh s})^{2} s^{2} ds d\tau).
\end{equation}
Choosing $R = \delta \frac{1}{\mathcal E(T)^{2}}$, $(\ref{4.31})$ can be absorbed into the left hand side of $(\ref{4.28.1})$, proving
\begin{equation}\label{4.28.2}
\aligned
\ln(\mathcal E(T)) - \ln(\mathcal E(0)) \lesssim \frac{\epsilon^{2}}{\delta} (\ln(\frac{1}{\delta}) + \ln(\mathcal E(T))) + \epsilon^{2}.
\endaligned
\end{equation}
This implies a uniform bound on $\mathcal E(T)$. Plugging the uniform bound on $\mathcal E(\tau)$ for all $\tau$ further implies a uniform bound on
\begin{equation}\label{4.28.3}
\int_{0}^{T} \int (\frac{s}{\sinh s})^{2} (\frac{\cosh s}{\sinh s}) \tilde{v}(\tau, s)^{4} s^{2} ds d\tau < \infty.
\end{equation}
This proves scattering, assuming Lemma $\ref{l4.1}$ is true. $\Box$\vspace{5mm}

\noindent \emph{Proof of Lemma $\ref{l4.1}$:} For $t > 1$,
\begin{equation}\label{4.40}
u(t) = S(t)(u_{0}, u_{1}) + \int_{0}^{t} S(t - t')(0, u^{3}) dt' = u_{l} + u_{nl}.
\end{equation}
First take the Duhamel term $u_{nl}$. Because the curve $t^{2} - r^{2} = 1$ has slope $\frac{dr}{dt} > 1$ everywhere,
\begin{equation}\label{4.42}
s \tilde{u}_{nl}(\tau, s)|_{\tau = 0} = \int_{1}^{e^{\tau} \cosh s} \int_{e^{\tau} \sinh s - e^{\tau} \cosh s + t}^{e^{\tau} \sinh s + e^{\tau} \cosh s - t} r u^{3}(t,r) dr dt.
\end{equation}
By direct computation,
\begin{equation}\label{4.43}
\aligned
\int_{0}^{\infty} (\partial_{\tau}(s \tilde{u}_{nl})|_{\tau = 0})^{2} ds \lesssim \int_{0}^{\infty} e^{2s} (\int_{1}^{\cosh s} (e^{s} - t) u^{3}(t, e^{s} - t) dt)^{2} ds \\ + \int_{0}^{\infty} e^{-2s} (\int_{1}^{\cosh s} (t - e^{-s}) u^{3}(t, t - e^{-s}) dt)^{2} ds.
\endaligned
\end{equation}
By H{\"o}lder's inequality, since $e^{s} - \cosh s \sim e^{s}$,
\begin{equation}\label{4.44}
\aligned
\int_{0}^{\infty} e^{2s} (\int_{1}^{\cosh s} (e^{s} - t) u^{3}(t, e^{s} - t) dt)^{2} ds \\ \lesssim \int_{0}^{\infty} \int_{1}^{\cosh s} e^{3s} (e^{s} - t)^{2} u^{6}(t, e^{s} - t) dt ds \\ 
\lesssim \int_{0}^{\infty} \int_{t^{2} - r^{2} \leq 1} u^{6}(t, r) r^{4} dt dr < \infty.
\endaligned
\end{equation}
The last inequality follows from global well-posedness of $u$, which implies $\| u \|_{L_{t,x}^{4}([1, 3] \times \mathbf{R}^{3})} < \infty$, $(\ref{4.7})$, Strichartz estimates, and the radial Sobolev embedding theorem, which implies
\begin{equation}\label{4.45}
\| |x|^{1/3} u \|_{L_{t,x}^{6}(\mathbf{R} \times \mathbf{R}^{3})} \lesssim \| |\nabla|^{1/6} u \|_{L_{t}^{6} L_{x}^{3}(\mathbf{R} \times \mathbf{R}^{3})}.
\end{equation}
Also by a change of variables and H{\"o}lder's inequality, since $(t - e^{-s}) \gtrsim 1$ for $s \geq 1$ and $t \geq 1$,
\begin{equation}\label{4.46}
\aligned
\int_{1}^{\infty} e^{-2s} (\int_{1}^{\cosh s} (t - e^{-s}) u^{3}(t, t - e^{-s}) dt)^{2} ds \\ \lesssim \int_{1}^{\infty} \int_{1}^{\cosh s} e^{-s} (t - e^{-s})^{2} u^{6}(t, e^{s} - t) dt ds \\ 
\lesssim \int_{0}^{\infty} \int_{t^{2} - r^{2} \leq 1} u^{6}(t, r) r^{4} dt dr < \infty.
\endaligned
\end{equation}
Also, by the radial Sobolev embedding theorem and Young's inequality,
\begin{equation}\label{4.47}
\aligned
\int_{1}^{\infty} e^{-2s} (\int_{1}^{\cosh s} (t - e^{-s}) u^{3}(t, t - e^{-s}) dt)^{2} ds \\ \lesssim \int_{1}^{3} (\int_{t^{2} - r^{2} \leq 1} u(t, r)^{6} r^{2} dr)^{1/2} dt \lesssim \int_{1}^{3} \frac{1}{(t - 1)^{3/4}} dt < \infty.
\endaligned
\end{equation}
This takes care of the nonlinear Duhamel piece.\vspace{5mm}

Now consider the linear piece. First consider the contribution of
\begin{equation}\label{4.48}
S(t - 1)(u_{0}, 0).
\end{equation}
Recall that if $w$ solves $(\ref{4.48})$ and $r > t$,
\begin{equation}\label{4.49}
r w(t,r) = \frac{1}{2} [u_{0}(t + r)(t + r) + u_{0}(r - t) (r - t)],
\end{equation}
so if $u_{1} = 0$, $u_{l} = S(t - 1)(u_{0}, 0)$, and
\begin{equation}\label{4.50}
\aligned
s \tilde{u}_{l}(\tau, s) = e^{\tau} \sinh s \cdot u_{l}(e^{\tau} \cosh s, e^{\tau} \sinh s) \\ = \frac{1}{2} [u_{0}(e^{\tau + s} - 1) (e^{\tau + s} - 1) + u_{0}(1 - e^{\tau - s}) (1 - e^{\tau - s})].
\endaligned
\end{equation}
Let $\chi \in C_{0}^{\infty}(\mathbf{R})$ be a function satisfying
\begin{equation}\label{4.51}
1 = \sum_{k \geq 0} \chi(s - k),
\end{equation}
for any $s \in [0, \infty)$, and $\chi(s - k)$ is supported on $(k - 1) \cdot \ln(2) \leq s \leq (k + 1) \cdot \ln(2)$, and split 
\begin{equation}
\tilde{u}_{l}(\tau, s) = \tilde{u}_{l}^{(1)}(\tau, s) + \tilde{u}_{l}^{(2)}(\tau, s) + \tilde{u}_{l}^{(3)}(\tau,s) + \tilde{u}_{l}^{(4)}(\tau, s) + \tilde{u}_{l}^{(5)}(\tau,s) + \tilde{u}_{l}^{(6)}(\tau, s),
\end{equation}
where
\begin{equation}\label{4.52}
\aligned
s \tilde{u}_{l}^{(1)}(\tau, s) = \sum_{k \geq 0} \chi(s - k) (P_{\leq -k} u_{0})(e^{\tau + s} - 1) \cdot (e^{\tau + s} - 1), \\
s \tilde{u}_{l}^{(2)}(\tau, s) = P_{\geq 0} \sum_{k \geq 0} \chi(s - k) (P_{> -k} u_{0})(e^{\tau + s} - 1) \cdot (e^{\tau + s} - 1), \\
s \tilde{u}_{l}^{(3)}(\tau, s) = P_{\leq 0} \sum_{k \geq 0} \chi(s - k) (P_{> -k} u_{0})(e^{\tau + s} - 1) \cdot (e^{\tau + s} - 1), \\
s \tilde{u}_{l}^{(4)}(\tau, s) = \sum_{k \geq 0} \chi(s - k) (P_{\leq k} u_{0})(1 - e^{\tau - s}) \cdot (1 - e^{\tau - s}), \\
s \tilde{u}_{l}^{(5)}(\tau, s) = P_{\leq 0} \sum_{k \geq 0} \chi(s - k) (P_{> k} u_{0})(1 - e^{\tau - s}) \cdot (1 - e^{\tau - s}), \\
s \tilde{u}_{l}^{(6)}(\tau, s) = P_{\geq 0} \sum_{k \geq 0} \chi(s - k) (P_{> k} u_{0})(1 - e^{\tau - s}) \cdot (1 - e^{\tau - s}).
\endaligned
\end{equation}
Taking the derivative,
\begin{equation}\label{4.53}
\aligned
\partial_{\tau} (s \tilde{u}_{l}^{(1)})(\tau, s)|_{\tau = 0} =  \sum_{k \geq 0} \chi(s - k) (P_{\leq -k} u_{0}')(e^{s} - 1) \cdot (e^{s} - 1) e^{s} \\
+ \sum_{k \geq 0} \chi(s - k) (P_{\leq -k} u_{0})(e^{s} - 1) \cdot e^{s}.
\endaligned
\end{equation}
Then by a change of variables, Hardy's inequality, and Young's inequality,
\begin{equation}\label{4.54}
\aligned
\| (\ref{4.53}) \|_{L^{2}[0, \infty)} \lesssim (\sum_{k \geq 0} 2^{k} (\sum_{j \leq -k} \| \chi(s - k)  (P_{j} \nabla u_{0})(e^{s} - 1) \|_{L^{2}} \\ + \| \chi(s - k) \frac{1}{|x|} (P_{j} u_{0})(e^{s} - 1) \|_{L^{2}})^{2})^{1/2} \lesssim \| u_{0} \|_{\dot{H}^{1/2}}.
\endaligned
\end{equation}
The computation of $\partial_{s}(s \tilde{u}_{l}^{(1)}(\tau, s))|_{\tau = 0}$ is similar, except that, in addition, it is necessary to compute
\begin{equation}\label{4.55}
\sum_{k} \| \chi'(s - k) (P_{\leq -k} u_{0})(e^{s} - 1) \cdot (e^{s} - 1) \|_{L^{2}}^{2}.
\end{equation}
Again, by a change of variables,
\begin{equation}\label{4.56}
(\ref{4.55}) \lesssim \sum_{k \geq 0} 2^{k} (\sum_{j \leq -k} \| \chi'(s - k) \frac{1}{|x|} (P_{j} u_{0})(e^{s} - 1) \|_{L^{2}}^{2} \lesssim \| u_{0} \|_{\dot{H}^{1/2}}^{2}.
\end{equation}
By the product rule,
\begin{equation}\label{4.57}
s \partial_{s} \tilde{u}_{l}(\tau, s) = \partial_{s}(s \tilde{u}_{l}(\tau, s)) - \tilde{u}_{l}(\tau, s).
\end{equation}
By the support properties of $\chi(s - k)$ and the Sobolev embedding theorem,
\begin{equation}\label{4.58}
\| \sum_{k \geq 0} \chi(s - k) (P_{\leq -k} u_{0})(e^{s} - 1) \cdot (e^{s} - 1) \|_{L^{\infty}} \lesssim \| u_{0} \|_{\dot{H}^{1/2}},
\end{equation}
and therefore,
\begin{equation}\label{4.59}
\| \frac{1}{s} \sum_{k \geq 2} \chi(s - k) (P_{\leq -k} u_{0})(e^{s} - 1) \cdot (e^{s} - 1) \|_{L^{2}([0, \infty)} \lesssim (\int_{1}^{\infty} \frac{1}{s^{2}} ds)^{1/2} \| u_{0} \|_{\dot{H}^{1/2}} \lesssim \| u_{0} \|_{\dot{H}^{1/2}}.
\end{equation}
Also, by the support properties of $\chi(s - k)$ and $(\ref{4.58})$,
\begin{equation}\label{4.60}
\| \sum_{k = 0,1} \chi(s - k) P_{\leq -k} u_{0}(e^{s} - 1) \cdot \frac{(e^{s} - 1)}{s} \|_{L^{2}([0, \infty)} \lesssim \| u_{0} \|_{\dot{H}^{1/2}}.
\end{equation}
Therefore, $\tilde{u}_{l}^{(1)}(\tau, s)|_{\tau = 0}$ has finite energy.\vspace{5mm}

Next, for any $k \geq 0$, $j > -k$, by the product rule and change of variables,
\begin{equation}\label{4.61}
\aligned
\| \partial_{\tau}(\chi(s - k) (P_{j} u_{0})(e^{s + \tau} - 1) \cdot (e^{s + \tau} - 1))|_{\tau = 0} \|_{L^{2}([0, \infty)} \\
\lesssim \| \chi(s - k) (P_{j} \nabla u_{0})(e^{s} - 1) \cdot (e^{s} - 1) e^{s} \|_{L^{2}([0, \infty)} \\
+ \| \chi(s - k) (P_{j} u_{0})(e^{s} - 1) \cdot e^{s} \|_{L^{2}([0, \infty)} \\
\lesssim 2^{k/2} \| P_{j} \nabla u_{0} \|_{L^{2}(2^{k - 1} - 1 \leq r \leq 2^{k + 1})} + 2^{k/2} \| \frac{1}{|x|} P_{j} u_{0} \|_{L^{2}(2^{k - 1} - 1 \leq r \leq 2^{k + 1})}.
\endaligned
\end{equation}
Therefore, if $f \in \dot{H}^{1/2}(\mathbf{R}^{3})$ is a radial function, by Bernstein's inequality,
\begin{equation}\label{4.62}
\aligned
\int_{0}^{\infty} (P_{l} f(s)) s \cdot \partial_{\tau}(\chi(s - k) (P_{j} u_{0})(e^{s + \tau} - 1) \cdot (e^{s + \tau} - 1))|_{\tau = 0} ds \\
\lesssim \| P_{l} f \|_{L^{2}(2^{k - 1} - 1 \leq r \leq 2^{k + 1})} [2^{k/2} \| P_{j} \nabla u_{0} \|_{L^{2}(2^{k - 1} - 1 \leq r \leq 2^{k + 1})} \\ + 2^{k/2} \| \frac{1}{|x|} P_{j} u_{0} \|_{L^{2}(2^{k - 1} - 1 \leq r \leq 2^{k + 1})}].
\endaligned
\end{equation}
Summing up, by Young's inequality, Bernstein's inequality,
\begin{equation}\label{4.63}
\aligned
\sum_{l \geq j + k > 0} \| P_{l} f \|_{L^{2}(2^{k - 1} - 1 \leq r \leq 2^{k + 1})} [2^{k/2} \| P_{j} \nabla u_{0} \|_{L^{2}(2^{k - 1} - 1 \leq r \leq 2^{k + 1})} \\ + 2^{k/2} \| \frac{1}{|x|} P_{j} u_{0} \|_{L^{2}(2^{k - 1} - 1 \leq r \leq 2^{k + 1})}] \lesssim \| f \|_{\dot{H}^{1/2}} \| u_{0} \|_{\dot{H}^{1/2}}.
\endaligned
\end{equation}
Next, by a change of variables,
\begin{equation}\label{4.64}
\| \chi(s - k) (P_{j} u_{0})(e^{s} - 1) \cdot (e^{s} - 1) \|_{L^{2}([0, \infty)} \lesssim 2^{-k/2} \| P_{j} u_{0} \|_{L^{2}(2^{k - 1} - 1 \leq r \leq 2^{k + 1})}.
\end{equation}
By the product rule,
\begin{equation}\label{4.65}
\aligned
 \partial_{\tau}(\chi(s - k) (P_{j} u_{0})(e^{s + \tau} - 1) \cdot (e^{s + \tau} - 1))|_{\tau = 0} \\ = \partial_{s} (\chi(s - k) (P_{j} u_{0})(e^{s} - 1) \cdot (e^{s} - 1)) - \chi'(s - k) (P_{j} u_{0})(e^{s} - 1) \cdot (e^{s} - 1).
 \endaligned
\end{equation}
Integrating by parts,
\begin{equation}\label{4.66}
\aligned
\int_{0}^{\infty} (P_{l} f(s)) s \cdot \partial_{s}(\chi(s - k) (P_{j} u_{0})(e^{s} - 1) \cdot (e^{s} - 1)) ds \\
= -\int_{0}^{\infty} [(P_{l} \nabla f(s)) s + (P_{l} f(s))] \chi(s - k) (P_{j} u_{0})(e^{s} - 1) \cdot (e^{s} - 1) ds \\
\lesssim 2^{-k/2} [\| P_{l} \nabla f \|_{L^{2}(2^{k - 1} - 1 \leq r \leq 2^{k + 1})} + \| \frac{1}{|x|} P_{l} f \|_{L^{2}(2^{k - 1} - 1 \leq r \leq 2^{k + 1})}] \| P_{j} u_{0} \|_{L^{2}(2^{k - 1} - 1 \leq r \leq 2^{k + 1})}.
\endaligned
\end{equation}
Summing up,
\begin{equation}\label{4.67}
\aligned
\sum_{0 \leq l < j + k} 2^{-k/2} [\| P_{l} \nabla f \|_{L^{2}(2^{k - 1} - 1 \leq r \leq 2^{k + 1})} \\ + \| P_{l} f \|_{L^{2}(2^{k - 1} - 1 \leq r \leq 2^{k + 1})}] \| P_{j} u_{0} \|_{L^{2}(2^{k - 1} - 1 \leq r \leq 2^{k + 1})} \lesssim \| f \|_{\dot{H}^{1/2}} \| u_{0} \|_{\dot{H}^{1/2}}.
\endaligned
\end{equation}
Also,
\begin{equation}\label{4.68}
\aligned
\int_{0}^{\infty} (P_{l} f(s))s \cdot \chi'(s - k) (P_{j} u_{0})(e^{s} - 1) \cdot (e^{s} - 1) ds \\ \lesssim \| P_{l} f \|_{L^{2}(2^{k - 1} - 1 \leq s \leq 2^{k + 1})} 2^{-k/2} \| P_{j} u_{0} \|_{L^{2}(2^{k - 1} - 1 \leq s \leq 2^{k + 1})}.
\endaligned
\end{equation}
Then by Bernstein's inequality,
\begin{equation}\label{4.69}
\aligned
\sum_{0 \leq l < j + k} 2^{-k/2} \| P_{l} f \|_{L^{2}(2^{k - 1} - 1 \leq r \leq 2^{k + 1}} \| P_{j} u_{0} \|_{L^{2}(2^{k - 1} - 1 \leq r \leq 2^{k + 1})} \lesssim \| f \|_{\dot{H}^{1/2}} \| u_{0} \|_{\dot{H}^{1/2}}.
\endaligned
\end{equation}
Therefore,
\begin{equation}\label{4.70}
\| \partial_{\tau}(\tilde{u}_{l}^{(2)}(\tau, s))|_{\tau = 0} \|_{\dot{H}^{-1/2}(\mathbf{R}^{3})} \lesssim \| u_{0} \|_{\dot{H}^{1/2}}.
\end{equation}
Also, by the product rule, and a change of variables
\begin{equation}\label{4.71}
\aligned
\| \partial_{s}(\chi(s - k) (P_{j} u_{0})(e^{s} - 1) \cdot (e^{s} - 1)) \|_{L^{2}} \lesssim 2^{-k/2} \| P_{j} u_{0} \|_{L^{2}(2^{k - 1} - 1 \leq r \leq 2^{k + 1})} \\
+ 2^{k/2} \| (P_{j} \nabla u_{0}) \|_{L^{2}(2^{k - 1} - 1 \leq r \leq 2^{k + 1})} + 2^{k/2} \| \frac{1}{|x|} (P_{j} u_{0}) \|_{L^{2}(2^{k - 1} - 1 \leq r \leq 2^{k + 1})}.
\endaligned
\end{equation}
Meanwhile,
\begin{equation}\label{4.71.1}
\| \frac{1}{s} \chi(s - k) (P_{j} u_{0})(e^{s} - 1) \cdot (e^{s} - 1) \|_{L^{2}} \lesssim 2^{-k/2} \| P_{j} u_{0} \|_{L^{2}(2^{k - 1} - 1 \leq r \leq 2^{k + 1})}.
\end{equation}
Then by Bernstein's inequality and Young's inequality,
\begin{equation}\label{4.72}
\aligned
\sum_{l} \| P_{l}(\sum_{l \leq k + j, k + j > 0} \chi(s - k) (P_{j} u_{0})(e^{s} - 1) \cdot (e^{s} - 1)) \|_{L^{2}}^{2} \\ \lesssim \sum_{l} 2^{l} \sum_{k} (\sum_{l \leq k + j, k + j > 0} \| \chi(s - k) (P_{j} u_{0})(e^{s} - 1) \cdot (e^{s} - 1) \|_{L^{2}})^{2} \lesssim \| u_{0} \|_{\dot{H}^{1/2}}^{2}.
\endaligned
\end{equation}
Also by Bernstein's inequality and $(\ref{4.71})$,
\begin{equation}\label{4.73}
\aligned
\sum_{l} \| P_{l}(\sum_{0 < k + j < l} \chi(s - k) (P_{j} u_{0})(e^{s} - 1) \cdot (e^{s} - 1)) \|_{L^{2}}^{2} \\ \lesssim \sum_{l} 2^{l} \sum_{k} (\sum_{0 < k + j < l} \| \chi(s - k) (P_{j} u_{0})(e^{s} - 1) \cdot (e^{s} - 1) \|_{L^{2}})^{2} \lesssim \| u_{0} \|_{\dot{H}^{1/2}}^{2}.
\endaligned
\end{equation}
Therefore, we have proved,
\begin{equation}\label{4.74}
\| \partial_{\tau}(\tilde{u}_{l}^{(2)}(\tau, s))|_{\tau = 0} \|_{\dot{H}^{-1/2}} + \| \tilde{u}_{l}^{(2)}(\tau, s)|_{\tau = 0} \|_{\dot{H}^{1/2}} \lesssim \| u_{0} \|_{\dot{H}^{1/2}}.
\end{equation}\vspace{5mm}

Next, following $(\ref{4.67})$--$(\ref{4.70})$ with $P_{l}$, $l \geq 0$ replaced by $P_{\leq 0}$ and $f \in L^{2}(\mathbf{R}^{3})$,
\begin{equation}\label{4.75}
\| \partial_{\tau}(\tilde{u}_{l}^{(3)}(\tau, s))|_{\tau = 0} \|_{L^{2}} + \| \tilde{u}_{l}^{(3)}(\tau, s)|_{\tau = 0} \|_{\dot{H}^{1}} \lesssim \| u_{0} \|_{\dot{H}^{1/2}}.
\end{equation}\vspace{5mm}

Next consider $\tilde{u}_{l}^{(4)}(\tau, s)$. By the product rule,
\begin{equation}\label{4.76}
\aligned
\partial_{\tau} (s \tilde{u}_{l}^{(4)}(\tau, s))|_{\tau = 0} = - \sum_{k \geq 0} \chi(s - k) (P_{\leq k} \nabla u_{0})(1 - e^{-s}) \cdot (1 - e^{-s}) e^{-s} \\
- \sum_{k \leq 0} \chi(s - k) (P_{\leq k} u_{0})(1 - e^{-s}) e^{-s}.
\endaligned
\end{equation}
Then, by Young's inequality,
\begin{equation}\label{4.77}
\aligned
\| \partial_{\tau}(s \tilde{u}_{l}^{(4)}(\tau, s))|_{\tau = 0} \|_{L^{2}([0, \infty)} \lesssim \sum_{k \geq 0} 2^{-k} (\sum_{j \leq k} \| \nabla P_{j} u_{0} \|_{L^{2}(1 - 2^{-k - 1} \leq r \leq 1 - 2^{-k + 1})})^{2} \\
+ \sum_{k \geq 0} 2^{-k} (\sum_{j \leq k} \| \frac{1}{|x|} P_{j} u_{0} \|_{L^{2}(1 - 2^{-k - 1} \leq r \leq 1 - 2^{-k + 1})})^{2} \lesssim \| u_{0} \|_{\dot{H}^{1/2}}^{2}.
\endaligned
\end{equation}
Also, by the product rule,
\begin{equation}\label{4.78}
\partial_{s}(s \tilde{u}_{l}^{(4)}(\tau, s)) = -\partial_{\tau}(s \tilde{u}_{l}^{(4)}(\tau, s)) + \sum_{k \geq 0} \chi'(s - k) (P_{k} u_{0})(1 - e^{-s}) \cdot (1 - e^{-s}).
\end{equation}
Then by the finite overlapping property of $\chi(s - k)$ and the radial Sobolev embedding theorem,
\begin{equation}\label{4.79}
\| \sum_{k \geq 0} \chi'(s - k) (P_{k} u_{0})(1 - e^{-s}) \cdot (1 - e^{-s}) \|_{L^{2}([0, \infty)}^{2} \lesssim \sum_{k \geq 0} \| P_{k} u_{0} \|_{\dot{H}^{1/2}}^{2} \lesssim \| u_{0} \|_{\dot{H}^{1/2}}^{2}.
\end{equation}
Therefore,
\begin{equation}\label{4.80}
\| \partial_{s}(s \tilde{u}_{l}^{(4)}(\tau, s))|_{\tau = 0} \|_{L^{2}([0, \infty)} + \| \partial_{\tau}(s \tilde{u}_{l}^{(4)}(\tau, s))|_{\tau = 0} \|_{L^{2}([0, \infty)} \lesssim \| u_{0} \|_{\dot{H}^{1/2}}.
\end{equation}

Next, by a change of variables,
\begin{equation}\label{4.81}
\| \chi(s - k) (P_{j} u_{0})(1 - e^{-s}) \cdot (1 - e^{-s}) \|_{L^{2}} \lesssim 2^{k/2} \| P_{j} u_{0} \|_{L^{2}(1 - 2^{-k - 1} \leq r \leq 1 - 2^{-k + 1})}.
\end{equation}
Therefore, by Young's inequality,
\begin{equation}\label{4.82}
\| s \tilde{u}_{l}^{5}(\tau, s)|_{\tau = 0} \|_{L^{2}([0, \infty)}^{2} \lesssim \sum_{k \geq 0} 2^{k} (\sum_{j > k} \| P_{j} u_{0} \|_{L^{2}(1 - 2^{-k - 1} \leq r \leq 1 - 2^{-k + 1})})^{2} \lesssim \| u_{0} \|_{\dot{H}^{1/2}}^{2}.
\end{equation}
Therefore, by the Fourier support of $\tilde{u}_{l}^{(5)}$,
\begin{equation}\label{4.83}
\| \tilde{u}_{l}^{(5)}(\tau, s)|_{\tau = 0} \|_{\dot{H}^{1}(\mathbf{R}^{3})} \lesssim \| u_{0} \|_{\dot{H}^{1/2}(\mathbf{R}^{3})}.
\end{equation}
Also, if $f \in L^{2}$ and $f$ is supported on $|\xi| \leq 1$,
\begin{equation}\label{4.84}
\aligned
\int_{0}^{\infty} f(s) s \cdot \partial_{\tau}(s \tilde{u}_{l}^{(5)}(\tau, s))|_{\tau = 0} ds = -\int_{0}^{\infty} f(s) s \cdot \partial_{s}(s \tilde{u}_{l}(\tau, s))|_{\tau = 0} ds \\
- \int_{0}^{\infty} f(s) s \cdot \sum_{k \geq 0} \chi'(s - k) (P_{\geq k} u_{0})(1 - e^{-s}) \cdot (1 - e^{-s}) ds.
\endaligned
\end{equation}
Integrating by parts, by $(\ref{4.82})$,
\begin{equation}\label{4.85}
-\int_{0}^{\infty} f(s) s \cdot \partial_{s}(s \tilde{u}_{l}^{(5)}(\tau, s))|_{\tau = 0} ds = \int_{0}^{\infty} \partial_{s}(f(s) s) \cdot s \tilde{u}_{l}^{(5)}(\tau, s)|_{\tau = 0} ds \lesssim \| f \|_{L^{2}} \| u_{0} \|_{\dot{H}^{1/2}}.
\end{equation}
Also, by $(\ref{4.82})$,
\begin{equation}\label{4.86}
\int_{0}^{\infty} f(s) s \cdot \sum_{k \geq 0} \chi'(s - k) (P_{\geq k} u_{0})(1 - e^{-s}) \cdot (1 - e^{-s}) ds \lesssim \| f \|_{L^{2}} \| u_{0} \|_{\dot{H}^{1/2}}.
\end{equation}
Therefore,
\begin{equation}\label{4.87}
\| \partial_{\tau}(s \tilde{u}_{l}^{(5)}(\tau, s))|_{\tau = 0} \|_{L^{2}([0, \infty)} + \| \partial_{s}(s \tilde{u}_{l}^{(5)}(\tau, s))|_{\tau = 0} \|_{L^{2}([0, \infty)} \lesssim \| u_{0} \|_{\dot{H}^{1/2}}.
\end{equation}\vspace{5mm}

Finally, take $\tilde{u}_{l}^{(6)}(\tau, s)$. Take $f \in \dot{H}^{1/2}$ supported in Fourier space on $|\xi| \geq 1$. Then by the product rule and $(\ref{4.82})$,
\begin{equation}\label{4.88}
\aligned
\| \partial_{s}(\chi(s - k) (P_{j} u_{0})(1 - e^{-s}) \cdot (1 - e^{-s})) \|_{L^{2}([0, \infty)} \\ \lesssim 2^{k/2} \| P_{j} u_{0} \|_{L^{2}(1 - 2^{-k - 1} \leq r \leq 1 - 2^{-k + 1})} \\
+ 2^{-k/2} \| P_{j} \nabla u_{0} \|_{L^{2}(1 - 2^{-k - 1} \leq r \leq 1 - 2^{-k + 1})} + 2^{-k/2} \| \frac{1}{|x|} P_{j} u_{0} \|_{L^{2}(1 - 2^{-k - 1} \leq r \leq 1 - 2^{-k + 1})}.
\endaligned
\end{equation}
Also, by $(\ref{4.82})$ and $(\ref{4.81})$,
\begin{equation}\label{4.89}
\| \frac{1}{s} \chi(s - k) (P_{j} u_{0})(1 - e^{-s}) \cdot (1 - e^{-s}) \|_{L^{2}([0, \infty)} \lesssim 2^{k/2} \| P_{j} u_{0} \|_{L^{2}(1 - 2^{-k - 1} \leq r \leq 1 - 2^{-k + 1})}.
\end{equation}
Therefore, by Young's inequality,
\begin{equation}\label{4.90}
\sum_{l < j + k} 2^{l} \sum_{k} (\sum_{j > k} \| \chi(s - k) (P_{j} u_{0})(1 - e^{-s}) \cdot (1 - e^{-s}) \|_{L^{2}})^{2} \lesssim \| u_{0} \|_{\dot{H}^{1/2}}^{2}.
\end{equation}
Also, by Bernstein's inequality,
\begin{equation}\label{4.91}
\sum_{l \geq j + k} 2^{-l} \sum_{k} (\sum_{j > k} \| \nabla \chi(s - k) (P_{j} u_{0})(1 - e^{-s}) \cdot (1 - e^{-s}) \|_{L^{2}})^{2} \lesssim \| u_{0} \|_{\dot{H}^{1/2}}^{2}.
\end{equation}
Therefore, we have finally proved that if $u_{1} = 0$,
\begin{equation}\label{4.92}
\tilde{u}_{l}(\tau, s)|_{\tau = 0} \in \dot{H}^{1/2}(\mathbf{R}^{3}) + \dot{H}^{1}(\mathbf{R}^{3}),
\end{equation}
and
\begin{equation}\label{4.93}
\partial_{\tau}(\tilde{u}_{l}(\tau, s))|_{\tau = 0} \in \dot{H}^{-1/2}(\mathbf{R}^{3}) + L^{2}(\mathbf{R}^{3}).
\end{equation}\medskip

To compute the contribution of
\begin{equation}\label{4.94}
S(t)(0, u_{1})
\end{equation}
to $\tilde{u}_{l}(\tau, s)$, observe that
\begin{equation}\label{4.95}
\frac{\sin(t \sqrt{-\Delta})}{\sqrt{-\Delta}} f = \partial_{t} (\frac{\cos(t \sqrt{-\Delta})}{\Delta} f).
\end{equation}
Plugging in the formula for a solution to the wave equation when $r > t$, let $w(t,r) = \cos(t \sqrt{-\Delta}) f$. Then,
\begin{equation}\label{4.96}
\aligned
\partial_{t} (w(t,r)) = \frac{1}{2r} \partial_{t} (f(t + r)(t + r) + f(r - t) (r - t)) \\ = \frac{1}{2r} [f(t + r) + f'(t + r)(t + r) - f(r - t) - f'(r - t) (r - t)].
\endaligned
\end{equation}
Then decompose $\tilde{u}_{l}(\tau, s) = \tilde{u}_{l}^{(1)}(\tau, s) + \tilde{u}_{l}^{(2)}(\tau, s) + \tilde{u}_{l}^{(3)}(\tau, s)$, where
\begin{equation}\label{4.96.1}
\aligned
s \tilde{u}_{l}^{(1)}(\tau, s) = \frac{1}{2} [f'(e^{\tau + s} - 1) \cdot (e^{\tau + s} - 1) - f'(1 - e^{\tau - s}) \cdot (1 - e^{\tau - s})], \\
s \tilde{u}_{l}^{(2)}(\tau, s) = \frac{1}{2} (1 - \chi(s)) [f(e^{\tau + s} - 1) - f(1 - e^{\tau - s})], \\
s \tilde{u}_{l}^{(3)}(\tau, s) = \frac{1}{2} \chi(s) [f(e^{\tau + s} - 1) - f(1 - e^{\tau - s})].
\endaligned
\end{equation}

Since
\begin{equation}\label{4.97}
f = \frac{u_{1}}{\Delta} \in \dot{H}^{3/2}(\mathbf{R}^{3}),
\end{equation}
the contribution of
\begin{equation}\label{4.98}
f'(e^{\tau + s} - 1) \cdot (e^{\tau + s} - 1) - f'(1 - e^{\tau - s}) \cdot (1 - e^{\tau - s}),
\end{equation}
to 
\begin{equation}\label{4.99}
(\tilde{u}_{l}(\tau, s)|_{\tau = 0}, \partial_{\tau} \tilde{u}_{l}(\tau, s)|_{\tau = 0})
\end{equation}
may be analyzed in exactly the same manner as the contribution of $S(t)(u_{1}, 0)$. Therefore,
\begin{equation}\label{4.99.1}
\tilde{u}_{l}^{(1)}(\tau, s)|_{\tau = 0} \in \dot{H}^{1/2} + \dot{H}^{1},
\end{equation}
and
\begin{equation}\label{4.99.2}
\partial_{\tau}(\tilde{u}_{l}^{(1)}(\tau, s))|_{\tau = 0} \in \dot{H}^{-1/2} + L^{2}.
\end{equation}\vspace{5mm}

Next take $\tilde{u}_{l}^{(2)}(\tau, s)$. By a change of variables,
\begin{equation}\label{4.100}
\int_{1}^{\infty} (\partial_{s} f(e^{s} - 1))^{2} ds = \int_{1}^{\infty} (f'(e^{s} - 1) \cdot e^{s})^{2} ds \lesssim \int |f'(r)|^{2} r dr \lesssim \| f \|_{\dot{H}^{3/2}(\mathbf{R}^{3})}^{2},
\end{equation}
and
\begin{equation}\label{4.101}
\int_{1}^{\infty} (\partial_{s} f(1 - e^{-s}))^{2} ds = \int_{1}^{\infty} (f'(1 - e^{-s}) \cdot e^{-s})^{2} ds \lesssim |f'(r)|^{2} r dr \lesssim \| f \|_{\dot{H}^{3/2}(\mathbf{R}^{3})}^{2}.
\end{equation}
By an identical calculation,
\begin{equation}\label{4.102}
\int_{1}^{\infty} (\partial_{\tau} f(e^{s + \tau} - 1)|_{\tau = 0})^{2} ds = \int_{1}^{\infty} (f'(e^{s} - 1) \cdot e^{s})^{2} ds \lesssim \int |f'(r)|^{2} r dr \lesssim \| f \|_{\dot{H}^{3/2}(\mathbf{R}^{3})}^{2},
\end{equation}
and
\begin{equation}\label{4.103}
\int_{1}^{\infty} (\partial_{s} f(1 - e^{\tau - s})|_{\tau = 0})^{2} ds = \int_{1}^{\infty} (f'(1 - e^{-s}) \cdot e^{-s})^{2} ds \lesssim \int |f'(r)|^{2} r dr \lesssim \| f \|_{\dot{H}^{3/2}(\mathbf{R}^{3})}^{2}.
\end{equation}
Next, by the fundamental theorem of calculus, for $s_{0} \sim 1$,
\begin{equation}\label{4.104}
s_{0} [f(e^{s_{0}} - 1) - f(1 - e^{-s_{0}})]^{2} = s_{0} [\int_{1 - e^{-s_{0}}}^{e^{s_{0}} - 1} f'(r) dr]^{2} \lesssim \int |f'(r)|^{2} r dr \lesssim \| f \|_{\dot{H}^{3/2}}^{2}.
\end{equation}
Therefore, by $(\ref{4.101})$ and $(\ref{4.102})$,
\begin{equation}\label{4.105}
\| \partial_{\tau}(\tilde{u}_{l}^{(2)}(\tau, s))|_{\tau = 0} \|_{L^{2}} \lesssim \| f \|_{\dot{H}^{3/2}},
\end{equation}
and
\begin{equation}\label{4.106}
\| \tilde{u}_{l}^{(2)}(0,s) \|_{\dot{H}^{1}} \lesssim \| f \|_{\dot{H}^{3/2}}.
\end{equation}\vspace{5mm}

Finally, consider
\begin{equation}\label{4.107}
f(e^{\tau + s} - 1) - f(1 - e^{\tau - s}),
\end{equation}
when $s < 1$. By direct computation,
\begin{equation}\label{4.108}
\partial_{\tau} [f(e^{\tau + s} - 1) - f(1 - e^{\tau - s})]|_{\tau = 0} = f'(e^{s} - 1) \cdot e^{s} + f'(1 - e^{-s}) \cdot e^{-s}.
\end{equation}
Then for $g \in \dot{H}^{1/2}$, by Hardy's inequality,
\begin{equation}\label{4.109}
\int f'(e^{s} - 1) \cdot e^{s} \cdot g(s) s ds + \int f'(1 - e^{-s}) \cdot e^{-s} \cdot g(s) s ds \lesssim \| f \|_{\dot{H}^{3/2}} \| g \|_{\dot{H}^{1/2}}.
\end{equation}
Also, by the fundamental theorem of calculus,
\begin{equation}\label{4.110}
\aligned
f(e^{s} - 1) - f(1 - e^{-s}) = \int_{s - \frac{s^{2}}{2} + \frac{s^{3}}{3!} - ...}^{s + \frac{s^{2}}{2} + \frac{s^{3}}{3!} + ...} f'(r) dr \\
= \int_{0}^{1} f'(s + \theta (\frac{s^{2}}{2} + \frac{s^{3}}{3!} + ...)) \cdot (\frac{s^{2}}{2} + \frac{s^{3}}{3!} + ...) d\theta \\ + \int_{-1}^{0} f'(s + \theta(\frac{s^{2}}{2} - \frac{s^{3}}{3!} + ...) \cdot (\frac{s^{2}}{2} + \frac{s^{3}}{3!} + ...) d\theta.
\endaligned
\end{equation}
Therefore, since $\chi(s)$ is supported on $s \leq 1$,
\begin{equation}\label{4.111}
\| f(e^{s} - 1) - f(1 - e^{-s}) \|_{\dot{H}^{1/2}} \lesssim \| f \|_{\dot{H}^{3/2}}.
\end{equation}
This proves that
\begin{equation}\label{4.112}
\| \tilde{u}_{l}^{(3)}(\tau, s)|_{\tau = 0} \|_{\dot{H}^{1/2}} + \| \partial_{\tau} \tilde{u}_{l}^{(3)}(\tau, s)|_{\tau = 0} \|_{\dot{H}^{-1/2}} \lesssim \| f \|_{\dot{H}^{3/2}}.
\end{equation}
This finally completes the proof of Lemma $\ref{l4.1}$. $\Box$

\section{Profile decomposition}
\emph{Proof of Theorem $\ref{t1.2}$:} This completes the proof that for any $(u_{0}, u_{1}) \in \dot{H}^{1/2} \times \dot{H}^{-1/2}$, $(\ref{1.1})$ has a global solution that scatters both forward and backward in time. To prove $(\ref{1.5})$, it remains to prove that for a sequence of initial data $(u_{n}^{0}, u_{n}^{1})$ and for any $A < \infty$,
\begin{equation}\label{5.1}
\| u_{0}^{n} \|_{\dot{H}^{1/2}} + \| u_{1}^{n} \|_{\dot{H}^{-1/2}} \leq A,
\end{equation}
\begin{equation}\label{5.2}
\| u^{n} \|_{L_{t,x}^{4}(\mathbf{R} \times \mathbf{R}^{3})} \leq f(A) < \infty,
\end{equation}
where $f : [0, \infty) \rightarrow [0, \infty)$, and $u^{n}$ is the solution to $(\ref{1.1})$ with initial data $(u_{0}^{n}, u_{1}^{n})$.

To prove this, make a profile decomposition.
\begin{theorem}[Profile decomposition]\label{t5.1}
Suppose that there is a uniformly bounded, radially symmetric sequence
\begin{equation}\label{5.3}
\| u_{0}^{n} \|_{\dot{H}^{1/2}(\mathbf{R}^{3})} + \| u_{1}^{n} \|_{\dot{H}^{-1/2}(\mathbf{R}^{3})} \leq A < \infty.
\end{equation}
Then there exists a subsequence, also denoted $(u_{0}^{n}, u_{1}^{n}) \subset \dot{H}^{1/2} \times \dot{H}^{-1/2}$ such that for any $N < \infty$,
\begin{equation}\label{5.4}
S(t)(u_{0}^{n}, u_{1}^{n}) = \sum_{j = 1}^{N} \Gamma_{n}^{j} S(t)(\phi_{0}^{j}, \phi_{1}^{j}) + S(t)(R_{0, n}^{N}, R_{1,n}^{N}),
\end{equation}
with
\begin{equation}\label{5.5}
\lim_{N \rightarrow \infty} \limsup_{n \rightarrow \infty} \| S(t)(R_{0,n}^{N}, R_{1,n}^{N}) \|_{L_{t,x}^{4}(\mathbf{R} \times \mathbf{R}^{3})} = 0.
\end{equation}
$\Gamma_{n}^{j} = (\lambda_{n}^{j}, t_{n}^{j})$ belongs to the group $(0, \infty) \times \mathbf{R}$, which acts by
\begin{equation}\label{5.6}
\Gamma_{n}^{j} F(t,x) = \lambda_{n}^{j} F(\lambda_{n}^{j} (t - t_{n}^{j}), \lambda_{n}^{j} x).
\end{equation}
The $\Gamma_{n}^{j}$ are pairwise orthogonal, that is, for every $j \neq k$,
\begin{equation}\label{5.7}
\lim_{n \rightarrow \infty} \frac{\lambda_{n}^{j}}{\lambda_{n}^{k}} + \frac{\lambda_{n}^{k}}{\lambda_{n}^{j}} + (\lambda_{n}^{j})^{1/2} (\lambda_{n}^{k})^{1/2} |t_{n}^{j} - t_{n}^{k}| = \infty.
\end{equation}
Furthermore, for every $N \geq 1$,
\begin{equation}\label{5.8}
\aligned
\| (u_{0, n}, u_{1, n}) \|_{\dot{H}^{1/2} \times \dot{H}^{-1/2}}^{2} = \sum_{j = 1}^{N} \| (\phi_{0}^{j}, \phi_{0}^{k}) \|_{\dot{H}^{1/2} \times \dot{H}^{-1/2}}^{2} \\ + \| (R_{0, n}^{N}, R_{1, n}^{N}) \|_{\dot{H}^{1/2} \times \dot{H}^{-1/2}}^{2} + o_{n}(1).
\endaligned
\end{equation}
\end{theorem}

Theorem $\ref{t5.1}$ gives the profile decomposition
\begin{equation}\label{5.9}
S(t)(u_{0}^{n}, u_{1}^{n}) = \sum_{j = 1}^{N} S(t - t_{n}^{j}) (\lambda_{n}^{j} \phi_{0}^{j}(\lambda_{n}^{j} x), (\lambda_{n}^{j})^{2} \phi_{1}^{j}(\lambda_{n}^{j} x)) + S(t)(R_{0, n}^{N}, R_{1,n}^{N}).
\end{equation}
In the course of proving Theorem $\ref{t5.1}$, \cite{Ramos} proved
\begin{equation}\label{5.10}
S(\lambda_{n}^{j} t_{n}^{j})(\frac{1}{\lambda_{n}^{j}} u_{0}^{n}(\frac{x}{\lambda_{n}^{j}}), \frac{1}{(\lambda_{n}^{j})^{2}} u_{1}^{n}(\frac{x}{\lambda_{n}^{j}})) \rightharpoonup \phi_{0}^{j}(x),
\end{equation}
weakly in $\dot{H}^{1/2}(\mathbf{R}^{3})$, and
\begin{equation}\label{5.11}
\partial_{t}S(t + \lambda_{n}^{j} t_{n}^{j})(\frac{1}{\lambda_{n}^{j}} u_{0}^{n}(\frac{x}{\lambda_{n}^{j}}), \frac{1}{(\lambda_{n}^{j})^{2}} u_{1}^{n}(\frac{x}{\lambda_{n}^{j}}))|_{t = 0} \rightharpoonup \phi_{1}^{j}(x)
\end{equation}
weakly in $\dot{H}^{-1/2}(\mathbf{R}^{3})$. Then after passing to a subsequence, $\lambda_{n}^{j} t_{n}^{j}$ converges to some $t^{j}$. Changing $(\phi_{0}^{j}, \phi_{1}^{j})$ to $S(-t^{j})(\phi_{0}^{j}, \phi_{1}^{j})$ and absorbing the error into $(R_{0, n}^{N}, R_{1, n}^{N})$,
\begin{equation}\label{5.10}
(\frac{1}{\lambda_{n}^{j}} u_{0}^{n}(\frac{x}{\lambda_{n}^{j}}), \frac{1}{(\lambda_{n}^{j})^{2}} u_{1}^{n}(\frac{x}{\lambda_{n}^{j}})) \rightharpoonup \phi_{0}^{j}(x),
\end{equation}
and
\begin{equation}\label{5.11}
\partial_{t}S(t)(\frac{1}{\lambda_{n}^{j}} u_{0}^{n}(\frac{x}{\lambda_{n}^{j}}), \frac{1}{(\lambda_{n}^{j})^{2}} u_{1}^{n}(\frac{x}{\lambda_{n}^{j}}))|_{t = 0} \rightharpoonup \phi_{1}^{j}(x).
\end{equation}
Then if $u^{j}$ is the solution to $(\ref{1.1})$ with initial data $(\phi_{0}^{j}, \phi_{1}^{j})$, then
\begin{equation}\label{5.12}
\| u^{j} \|_{L_{t,x}^{4}(\mathbf{R} \times \mathbf{R}^{3})} \leq M_{j}.
\end{equation}

Next, suppose that after passing to a subsequence, $\lambda_{n}^{j} t_{n}^{j} \nearrow +\infty$. Theorem $\ref{t4.0}$ also implies that for any $(\phi_{0}, \phi_{1}) \in \dot{H}^{1/2} \times \dot{H}^{-1/2}$, there exists a solution $u$ to $(\ref{1.1})$ that is globally well-posed and scattering, and furthermore, that $u$ scatters to $S(t)(\phi_{0}, \phi_{1})$ as $t \searrow -\infty$.
\begin{equation}\label{5.13}
\lim_{t \rightarrow -\infty} \| u - S(t)(\phi_{0}, \phi_{1}) \|_{\dot{H}^{1/2} \times \dot{H}^{-1/2}} = 0.
\end{equation}
Indeed, by Strichartz estimates, the dominated convergence theorem, and small data arguments, for some $T < \infty$ sufficiently large, $(\ref{1.1})$ has a solution $u$ on $(-\infty, -T]$ such that
\begin{equation}\label{5.14}
\| u \|_{L_{t,x}^{4}((-\infty, -T] \times \mathbf{R}^{3})} \lesssim \epsilon, \qquad (u(-T, x), u_{t}(-T, x)) = S(-T)(\phi_{0}, \phi_{1}).
\end{equation}
and by Strichartz estimates,
\begin{equation}\label{5.15}
\lim_{t \rightarrow +\infty} \| S(t)(u(-t), u_{t}(-t)) - (\phi_{0}, \phi_{1}) \|_{\dot{H}^{1/2} \times \dot{H}^{-1/2}} \lesssim \epsilon^{3}.
\end{equation}
Then by the inverse function theorem, there exists some $(u_{0}(-T), u_{1}(-T))$ such that $(\ref{1.1})$ has a solution that scatters backward in time to $S(t)(\phi_{0}, \phi_{1})$, and by Theorem $\ref{t4.0}$, this solution must also scatter forward in time. Therefore,
\begin{equation}\label{5.16}
S(-t_{n}^{j})(\lambda_{n}^{j} \phi_{0}^{j}(\lambda_{n}^{j} x), (\lambda_{n}^{j})^{2} \phi_{1}^{j}(\lambda_{n}^{j} x))
\end{equation}
converges strongly to
\begin{equation}\label{5.17}
(\lambda_{n}^{j} u^{j}(-\lambda_{n}^{j} t_{n}^{j}, \lambda_{n}^{j} x), (\lambda_{n}^{j})^{2} u_{t}^{j}(-\lambda_{n}^{j} t_{n}^{j}, \lambda_{n}^{j} x))
\end{equation}
in $\dot{H}^{1/2} \times \dot{H}^{-1/2}$, where $u^{j}$ is the solution to $(\ref{1.1})$ that scatters backward in time to $S(t)(\phi_{0}^{j}, \phi_{1}^{j})$, and the remainder may be absorbed into $(R_{0, n}^{N}, R_{1, n}^{N})$. In this case as well,
\begin{equation}\label{5.18}
\| u^{j} \|_{L_{t,x}^{4}(\mathbf{R} \times \mathbf{R}^{3})} \leq M_{j} < \infty.
\end{equation}
The proof for $\lambda_{n}^{j} t_{n}^{j} \searrow -\infty$ is similar.

Also, by $(\ref{5.8})$, there are only finitely many $j$ such that $\| \phi_{0}^{j} \|_{\dot{H}^{1/2}} + \| \phi_{1}^{j} \|_{\dot{H}^{-1/2}} > \epsilon$. For all other $j$, small data arguments imply
\begin{equation}\label{5.19}
\| u^{j} \|_{L_{t,x}^{4}(\mathbf{R} \times \mathbf{R}^{3})} \lesssim \| \phi_{0}^{j} \|_{\dot{H}^{1/2}} + \| \phi_{1}^{j} \|_{\dot{H}^{-1/2}}.
\end{equation}
Then by the decoupling property $(\ref{5.7})$, $(\ref{5.12})$, $(\ref{5.19})$, and Lemma $\ref{l6.2}$,
\begin{equation}
\limsup_{n \nearrow \infty} \| u^{n} \|_{L_{t,x}^{4}(\mathbf{R} \times \mathbf{R}^{3})}^{2} \lesssim \sum_{j} \| u^{j} \|_{L_{t,x}^{4}(\mathbf{R} \times \mathbf{R}^{3})}^{2} < \infty.
\end{equation}
This proves Theorem $\ref{t1.2}$. $\Box$

\bibliographystyle{plain}

\end{document}